\documentclass[journal]{IEEEtran}  % Comment this line out

\usepackage{amsmath}
\usepackage{amssymb}
\usepackage{algorithm}
\usepackage{subcaption}
\usepackage[noend]{algpseudocode}
\usepackage[shortlabels]{enumitem}
\usepackage[colorinlistoftodos,prependcaption]{todonotes}
\IEEEoverridecommandlockouts                              % This command is only
\newcommand{\affop}{H^*}
\newcommand{\affconst}{T}
\algrenewcommand\textproc{\textsc}

\title{\LARGE \bf
A Unifying Complexity Certification Framework for Active-Set Methods for Convex Quadratic Programming}
\author{Daniel Arnstr\"om and Daniel Axehill% <-this % stops a space
  \thanks{This work was supported by the Swedish Research Council (VR) under contract number 2017-04710.}% <-this % stops a space
\thanks{D. Arnstr\"om and D. Axehill are with the Division of Automatic Control, Link\"oping University,
  Sweden 
{\tt\small daniel.\{arnstrom,axehill\}@liu.se}}%
}

\begin{document}
\renewcommand{\baselinestretch}{1.0}

\maketitle
\thispagestyle{empty}
\pagestyle{empty}
\newtheorem{proposition}{Proposition}
\newtheorem{lemma}{Lemma}
\newtheorem{corollary}{Corollary}
\newtheorem{remark}{Remark}
\newtheorem{theorem}{Theorem}
\newcommand\inv[1]{#1\raisebox{1.15ex}{$\scriptscriptstyle-\!1$}}
\newcommand\qrt[1]{#1\raisebox{1.15ex}{$\scriptscriptstyle-\!1/2$}}

\begin{abstract}
  In model predictive control (MPC) an optimization problem has to be solved at each time step, which in real-time applications makes it important to solve these optimization problems efficiently and to have good upper bounds on worst-case solution time. Often for linear MPC problems, the optimization problem in question is a quadratic program (QP) that depends on parameters such as system states and reference signals. A popular class of methods for solving such QPs is active-set methods, where a sequence of linear systems of equations is solved. We propose an algorithm for computing which sequence of subproblems an active-set algorithm will solve, for every parameter of interest. By knowing these sequences, a worst-case bound on how many iterations, and ultimately the maximum time, the active-set algorithm requires to converge can be determined. The usefulness of the proposed method is illustrated on a set of QPs, originating from MPC problems, by computing the exact worst-case number of iterations primal and dual active-set algorithms require to reach optimality.
\end{abstract}

\section{INTRODUCTION}
In model predictive control (MPC) an optimization problem has to be solved at each time step, which for linear MPC often is a quadratic program (QP) which depends on parameters such as system states and reference signals, making it a multi-parametric QP (mpQP). Often, these mpQPs are solved offline parametrically for a set of parameters and the pre-computed solution is then used online \cite{bemporad2002explicit}. However, the pre-computed solution grows exponentially in complexity with the dimensions of the problem and, for high-dimensional problems, limited memory can restrict the use of a pre-computed solution online. For such problems, the QP has to be solved online and the limited time and computational resources often at hand in real-time MPC require the employed QP solver to be efficient and to have guarantees on the time needed to solve the QPs within a given tolerance.

Popular methods for solving QPs encountered in MPC are active-set methods \cite{nocedal}\cite{10.1007BF02591962}\cite{10.1137S1052623400376135}\cite{qpoases}\cite{bemporad2015quadratic}, interior-point methods \cite{rao1998application}\cite{wang2010fast} and gradient projection methods \cite{10.1109TAC.2013.2275667}\cite{axehill2008dual}\cite{richter2009real}. Active-set methods easily integrate warm-starting of the solver, i.e., the use of a previous solution to start the solver in the next iteration, which often reduces the number of iterations needed by the solver \cite{10.1109TAC.2011.2108450}\cite{herceg2015dominant}. A well-known drawback of active-set methods is, however, that the complexity can be exponential in the worst-case \cite{klee1972good}. Although, polynomial complexity is often observed in practice \cite{spielman2004smoothed}. In contrast to active-set methods, theoretical polynomial bounds on the computational complexity of some interior-point and gradient projection methods have been proven in, e.g.,  \cite{rao1998application}\cite{10.1109TAC.2013.2275667}\cite{richter2012computational}\cite{giselsson2012execution}.

To close the gap between the possible exponential complexity and the often experienced polynomial complexity, methods for determining the exact complexity of the active-set QP methods presented in \cite{nocedal},\cite{10.1007BF02591962} and \cite{10.1137S1052623400376135} have been proposed in \cite{cdc2019},\cite{10.1109TAC.2017.2696742} and \cite{cimini2019complexity}, respectively. Similarly, a method for determining the complexity of a primal active-set methods for linear programs (LPs) has been proposed in \cite{10.1109TAC.2011.2108450}. 
This paper extends the result in \cite{cdc2019}, which handles the strictly convex case, to also handle positive semi-definite mpQPs, leading to additional theoretical as well as numerical results. In addition to being able to certify the complexity of primal active-set methods applied to positive semi-definite mpQPs, it is shown that this extension allows for dual active-set QP methods and active-set methods for linear programs to be certified with the presented method, enabling the results in \cite{cdc2019},\cite{10.1109TAC.2017.2696742} and \cite{10.1109TAC.2011.2108450} to be viewed in a unified framework. 

The main contribution of this paper is, hence, a method for analyzing \textit{exactly} which subproblems, i.e., systems of linear equations, a primal active-set algorithm will solve in order to compute an optimal solution for any set of parameters in an mpQP, which can ultimately be used to determine the worst-case computational complexity of the algorithm. The proposed method is used offline on a given mpQP, giving a priori knowledge about how the active-set algorithm will act when employed online such as a worst-case bound on the number of iterations. Furthermore, exact knowledge about the subproblems that can be encountered can be used to tailor the solver for the specific mpQP at hand. 

A challenging aspect of the analysis of the primal active-set QP algorithm considered in this work is that it turns out that all iterates are not necessarily affine in the parameter, in contrast to the methods studied in \cite{10.1109TAC.2017.2696742}, \cite{cimini2019complexity} and \cite{10.1109TAC.2011.2108450}. Nonaffine iterates are shown to lead to a partition of the parameter space consisting of both linear and quadratic inequalities, in contrast to only linear inequalities which is the case in \cite{10.1109TAC.2017.2696742}, \cite{cimini2019complexity} and \cite{10.1109TAC.2011.2108450}. 

The rest of the paper is outlined as follows: Section \ref{sec:prelim} introduces notation, some background theory and the active-set algorithm considered. Properties of this algorithm are then presented in Section \ref{sec:prop}  which are used in the proposed complexity certification method presented in Section \ref{sec:cert-primal-as}. The active-set algorithm as well as the certification method is extended in Section \ref{sec:singular} to also handle positive semi-definite QPs and how these extensions unify results for primal and dual active-set QP algorithms as well as active-set LP algorithms is discussed. Finally, in Section \ref{sec:examples} the proposed method is illustrated on a set of examples, including MPC problems that are representative of problems encountered in real-time MPC.
\newpage
\section{Preliminaries}
\label{sec:prelim}
It is well-known that a linear MPC problem can be cast into an mpQP on the form \eqref{eq:mpQP}, where the parameter $\theta$ contains the measured/estimated state \cite{bemporad2002explicit},  
\begin{equation}
  \label{eq:mpQP}
  \begin{aligned}
	&\underset{x}{\text{minimize}}&&\frac{1}{2}x^T H x + (f^T + \theta^T f_{\theta}^T) x \\
	&\text{subject to} &&Ax \leq b + W \theta.\\
  \end{aligned}
\end{equation}
Where $x \in \mathbb{R}^n$ and the parameter $\theta \in \Theta_0 \subseteq \mathbb{R}^p $, with $\Theta_0$ being a polyhedron. The mpQP is given by $A\in \mathbb{R}^{m\times n}$, $b \in \mathbb{R}^{m}$, $W\in \mathbb{R}^{m\times p}$, $f\in \mathbb{R}^{n}$, $f_{\theta}\in\mathbb{R}^{n\times p}$, and $H \in \mathbb{S}^{n}_{+}$.
For convenience, we also introduce the compact notation $b(\theta) = b+W\theta$ and $f(\theta) = f+f_{\theta} \theta$ which will sometimes be used to clean up expressions. 

The feasible set can also be expressed in terms of each constraint as $[A]_i x \leq [b]_i+[W]_i \theta, i \in \mathcal{K}$, where the notation $[.]_i$ means the $i$:th row of the corresponding matrix and $ {\mathcal{K}\triangleq\{1,2,...,m\}}$.   
A constraint is said to be \textit{active} if it holds with equality.

The primal active-set algorithm to be studied is an iterative algorithm which searches for the active constraints at the optimum, motivating the following notation. $x_k$ is the iterate at iteration $k$ and $\mathcal{W}_k$ is a subset of the constraints, called the working set, that are active at $x_k$. Moreover, we define $A_k, b_k$ and $W_k$ to denote the rows of the matrices indexed by $\mathcal{W}_k$ and we denote the complement of $\mathcal{W}_k$ as $\bar{\mathcal{W}}_k\triangleq \mathcal{K} \setminus \mathcal{W}_k$. The constrained set ${P_k \triangleq \{x\in \mathbb{R}^n|A_k x = b_k(\theta)\}}$ denotes the manifold defined by the working set at iteration $k$.

\subsection{Equality constrained mpQP}
The active-set algorithm considered in this paper solves a sequence of equality constrained QPs (EQPs) on the form 
\begin{equation}
  \label{eq:eqp}
  \begin{aligned}
	 &\underset{x}{\text{minimize}}&& x^T H x + f(\theta)^Tx\\
	 &\text{subject to} && A_k x =b_k(\theta) 
  \end{aligned}
\end{equation}
The optimizer $x^*_k$ of this subproblem, which we will call a \textit{constrained stationary point} (CSP), and the dual variable $\lambda_k$ can be obtained by solving the following linear system of equations, also known as a KKT-system, 
\begin{equation}
  \label{eq:kkt}
  \begin{pmatrix}
	H & A^T_k \\ A_k & 0 
  \end{pmatrix}
  \begin{pmatrix}
   x^*_k \\ \lambda_k 
\end{pmatrix} = 
\begin{pmatrix}
  -f(\theta)\\
  b_k(\theta)
\end{pmatrix}.
\end{equation}
If there exists a unique solution to \eqref{eq:kkt} the inverse of the KKT matrix can be partitioned as  
\begin{equation}
  \label{eq:kkt-inv}
  \begin{pmatrix}
	H & A^T_k \\ A_k & 0 
  \end{pmatrix}^{-1}  = 
  \begin{pmatrix}
	{H}^*_k & T_k \\ T_k^T & U_k  
  \end{pmatrix}
\end{equation}
and the solution to \eqref{eq:kkt} is given by 
\begin{equation}
  \label{eq:const-kkt}
  \begin{aligned}
	x^*_k &= -{H^*_k} f(\theta) + T_k b_k(\theta) \\
	\lambda_k &= -T^T_k f(\theta) + U_k b_k(\theta)
  \end{aligned}
\end{equation}
Importantly, the solution to the KKT-system in \eqref{eq:const-kkt} is affine in $\theta$, i.e., 
\begin{equation}
  \begin{aligned}
    \label{eq:kkt_para_sol}
	x^*_k = F^*_k \theta + G^*_k, \quad
	\lambda_k = F^{\lambda}_k \theta + G^{\lambda}_k 
  \end{aligned}
\end{equation}
with $F^*_k, G^*_k, F^{\lambda}_k ,G^{\lambda}_k$ defined by
\begin{subequations}
	\label{eq:kkt-explicit}
  \begin{align}
	\label{eq:F-star}F^*_k &\triangleq -H^*_k f_{\theta}+T_k W_k, &&G^*_k \triangleq -H^*_k f+T_k b_k \\
	\label{eq:F-lam}F^{\lambda}_k &\triangleq-T_k^T f_{\theta}+U_k W_k, &&G^{\lambda}_k \triangleq -T^T_k f + U_k b_k 
  \end{align}
\end{subequations}

When $H \succ 0$ and $A_k$ has full row rank, $H^*_k, T_k$ and $U_k$ can be expressed explicitly as  \cite{fletcher}
\begin{equation}
  \begin{aligned}
	H^*_k &= \inv{H} - \inv{H} A^T_k \inv{(A_k \inv{H} A^T_k)} A_k \inv{H}  \\  
	T_k &= \inv{H} A_k^T \inv{(A_k \inv{H} A_k^T)}\\  
	U_k &= -(A_k^T \inv{H} A^T_k) \\  
  \end{aligned}
\end{equation}
This representation is used in so-called range-space methods for solving the KKT-system. Evidently a range-space method can not be used when $H$ is singular since $\inv{H}$ is needed. Nevertheless, the KKT-matrix might still be non-singular if $H$ is positive definite on the nullspace of $A_k$. Formally this can be expressed as the \textit{reduced Hessian}  $Z_k^T H Z_k$ being positive definite, where $Z_k$ is a matrix with columns forming a basis for the nullspace of $A_k$. By introducing $Y_k$ as a matrix with columns spanning the range-space of $A_k$ and satisfying $Y_k^T A_k = I$, $H^*_k, T_k$ and $U_k$ can be expressed explicitly as \cite{fletcher} 

\begin{equation}
  \label{eq:Hstar-def}
  \begin{aligned}
	H^*_k &= Z_k\inv{(Z_k^T H Z_k)} Z_k^T  \\
	T_k &= Y_k-Z_k\inv{(Z_k^T H Z_k)} Z_k^T H Y_k \\ 
	U_k &= Y_k^T H Z_k\inv{(Z_k^T H Z_k)} Z_k^T H Y_k - Y_k^T H Y_k 
  \end{aligned}
\end{equation}
This representation is used in so-called null-space methods for solving the KKT-system.
\begin{remark}
  $Z_k^T H Z_k \succ 0$ is sufficient for the KKT-system \eqref{eq:kkt} to have a unique solution. In particular, note that 
  $H \succ 0$ and $A_k$ full row rank $\implies$  $Z^T_k H Z_k \succ 0$. 
\end{remark}

Since null-space methods encapsulate the semi-definite case, the formulations in \eqref{eq:Hstar-def} will be considered in the sequal. For the strictly convex case, however, all results can be translated to the case when a range-space method is used.  

Before proceeding, we prove the following projective property of $H^*$ which will be central when the properties of the active-set algorithm, soon to be introduced, is discussed in Section \ref{sec:prop}.
\begin{lemma}
  \label{lem:proj-operator}
  $P_{k+1} \subseteq P_k  \implies H^*_{k+1} H H^*_k = H^*_{k+1}$
\end{lemma}
\begin{IEEEproof}
  Since $P_{k+1}  \subseteq P_k$ we have that ${Z_{k} = [Z_{k+1},\: Z_+]}$. Using this together with the formula for the inverse of a 2x2 block matrix gives 
\begin{equation}
  \begin{aligned}
	\inv{(Z_k^T H Z_k)}
	&= \begin{bmatrix}
	  Z_{k\text{+}1}^T H Z_{k\text{+}1} & Z_{k\text{+}1}^T H Z_\text{+}\\
	  Z^T_\text{+} H Z_{k\text{+}1} & Z^T_\text{+} H Z_\text{+}
	\end{bmatrix}^{\text{-}1} \triangleq 
 \begin{bmatrix}
   U & V \\ V^T & W 
 \end{bmatrix}^{\text{-}1}\\ 
	&= \begin{bmatrix}
	  U^{-1}(I + V \tilde{S} V^T U^{-1}) & - U^{-1} V \tilde{S} \\
	  -\tilde{S} V^T U^{-1}&  \tilde{S}  
	\end{bmatrix}
  \end{aligned}
\end{equation}
with $\tilde{S} \triangleq \inv{(W-V^T U V)}$ being the inverse of a Schur complement. Multiplication with $Z_k^T$ from the right then gives
\begin{equation}
  \label{eq:block1}
  \inv{(Z_k^T H Z_k)} Z_k^T  =  \begin{bmatrix}
	\inv{U} \left(Z_{k\text{+}1}^T \text{+} V \tilde{S} \tilde{Z} \right) \\
	-\tilde{S} \tilde{Z}
  \end{bmatrix}
\end{equation}
with $\tilde{Z} \triangleq (V^T\inv{U} Z^T_{k+1}-Z_{+})$. By definition we have that 
\begin{equation}
  \label{eq:uv}
  Z^T_{k+1} H Z_k = [U,\: V].
\end{equation}
Hence, multiplying \eqref{eq:block1} from the left with \eqref{eq:uv} and recalling the definition of $H^*_k$ from \eqref{eq:Hstar-def} gives, after some cancellations, 
\begin{equation}
  \label{eq:key}
  Z^T_{k+1} H H^*_k = Z^T_{k+1} H Z_k (Z^T_k H Z_k)^{-1} Z_k^T = Z^T_{k+1} 
\end{equation}
Finally, we get the desired result by recalling the definition of $H^*_{k+1}$ from \eqref{eq:Hstar-def} and using \eqref{eq:key}
\begin{equation}
  \begin{aligned}
	H^*_{k+1} H H^*_k =&  Z^T_{k+1} (Z_{k+1}^T H Z_{k+1})^{-1} Z^T_{k+1} H H^*_k  \\
	=& Z_{k+1}^T(Z^T_{k+1} H Z_{k+1})^{-1} Z^T_{k+1}  \\
	=& H^*_{k+1}
  \end{aligned}
\end{equation}
\end{IEEEproof}%

\subsection{A primal active-set algorithm}
\label{sec:primal-as}
%\todo[inline]{Be consistent with the use of $b$ or $b(\theta)$ and $f$ or $f(\theta)$}
An important class of methods for solving QPs are active-set methods, which solve the QP by solving  a sequence of EQPs, i.e., systems of linear equations. There are plenty of different primal active-set methods in the litterature, e.g., \cite{nocedal}\cite{fletcher1971general}\cite{dantzig1963linear}, and numerous of these are equivalent \cite{best1984equivalence} in the sense that they produce the same iterates given the same starting conditions. In this paper we consider the primal active-set algorithm given by Algorithm \ref{alg:primal-as}, described in detail below. This algorithm formulation is chosen to make the certification method, described in Section \ref{sec:cert-primal-as}, more succinct and the definition of an iteration of the algorithm sound. 
However, it would be possible to instead consider any other equivalent formulation, such as any of the primal active-set methods cited above. For example, this is done in \cite{cdc2019} where the algorithm formulation presented in \cite[Sec.16-5]{nocedal} is considered. 

Algorithm \ref{alg:primal-as} works for strictly convex QPs and can be extended to work for convex QPs. However, we will start by considering the strictly convex case to ease the initial analysis and then extend it to the semi-definite case in Section \ref{sec:singular}.  
\begin{algorithm}[H]
  \caption{Primal Active-Set Method for QP}\label{alg:primal-as}
  \begin{algorithmic}[1]
	\Require $x_0, \mathcal{W}_0, k=1$, dual tolerance $\epsilon_d \geq 0$
	\Ensure $x_k^*, \lambda_k, \mathcal{W}_k $
	\State $s_0 \leftarrow b+W\theta-A x_0$  
	\While {true}
	%\If{$Z^T_k H Z_k$ singular} \textsc{SingularIteration()}
	%\Else
	\algrenewcommand{\alglinenumber}[1]{\color{red}\footnotesize#1:}
	\State Compute $p_k$ by solving \eqref{eq:kkt};\quad  
  \State $[\sigma_k]_{\mathcal{C}_k} \leftarrow [A]_{\mathcal{C}_k} p_k$%, \: [\sigma_k]_{\mathcal{W}_k} \leftarrow 0
  , \: $s^*_k \leftarrow s_k - \sigma_k$%,\: $[s^*_k]_{\mathcal{W}_k} \leftarrow 0$   
	%\If{$s^*_k$ \geq 0$} Compute $\lambda_k$ by solving \eqref{eq:dual_var} 
	\If{$s^*_k \geq 0$} 
	\algrenewcommand{\alglinenumber}[1]{\color{green!60!black}\footnotesize#1:}
	\State Compute $\lambda_k$ by solving \eqref{eq:dual_var} 
  	\If{$\lambda_k \geq -\epsilon_d$}
	\Return $x_k^*, \lambda_k , \mathcal{W}_k$  
	\Else \quad $l \leftarrow \underset{i\in \mathcal{W}_k }{\text{argmin }} [\lambda_k]_i$;
	\: $\mathcal{W}_{k+1} \leftarrow \mathcal{W}_k\setminus\{l\}$ \label{line:const-rm}
	\State \quad $x_{k+1} \leftarrow  x_k+p_k$; \quad $s_{k+1} \leftarrow s_k^*$
  \EndIf
	\algrenewcommand{\alglinenumber}[1]{\color{red}\footnotesize#1:}
	\Else\quad $m \leftarrow \underset{i \in \bar{\mathcal{W}}_k:[s^*_k]_i < 0}{\text{argmin }} \frac{[s_k]_i}{[\sigma_k]_i};\:\:\: \mathcal{W}_{k+1}  \leftarrow \mathcal{W}_k \cup \{m\}$  
\State \quad $x_{k+1} \leftarrow x_k + \alpha^m_k p_k$;\:\:\:\:\: $s_{k+1} \leftarrow s_k - \alpha^m_k \sigma_k$
\EndIf
	\algrenewcommand{\alglinenumber}[1]{\color{black}\footnotesize#1:}
	\State $k \leftarrow k +1$
	%\EndIf
  \EndWhile
  \end{algorithmic}
\end{algorithm}

Algorithm \ref{alg:primal-as} starts with a feasible point $x_0$ and a corresponding working set $\mathcal{W}_0$, containing a subset of the constraints that are active at $x_0$. 
\begin{remark}
We allow $x_0$ to be affine in the parameter $\theta$, i.e., $x_0 = F_0 \theta + G_0$.
\end{remark} 

In an iteration of the algorithm, constraints are added to or removed from the working set while maintaining primal feasibility and updating the iterate. The iterate is updated in a line search fashion, i.e., $x_{k+1} = x_k + \alpha_k p_k$, with the search direction $p_k$ and the step length $\alpha_k$, defined below. 

  The search direction $p_k$ is the Newton step direction given by $p_k \triangleq x_k^* -x_k$ where $x^*_k$ is the solution to the EQP in \eqref{eq:eqp}. Instead of solving \eqref{eq:eqp} to obtain $p_k$,  one can reformulate its KKT-system \eqref{eq:kkt} in terms of $p_k$ instead of $x^*_k$ according to 

\begin{equation}
  \label{eq:kkt-step}
  \begin{pmatrix}
	H & A^T_k \\ A_k & 0 
  \end{pmatrix}
  \begin{pmatrix}
   p_k \\ \lambda_k 
\end{pmatrix} = 
\begin{pmatrix}
  -(Hx_k +f(\theta))\\
  0
\end{pmatrix}
\end{equation}
to obtain $p_k$ directly. In an iteration, we want to retain primal feasibility in the iterate while trying to move along a line segment from $x_k$ to $x^*_k$. Such a move can be done if $x^*_k$ is primal feasible, i.e., if ${b(\theta)-Ax^*_k \geq 0}$. The following notations prove useful when talking about primal feasibility of $x^*_k$ 
\begin{equation}
  \label{eq:opt_slack}
  s^*_k \triangleq b(\theta)-Ax^*_k = b(\theta)-A x_k - A p_k = s_k - \sigma_k
\end{equation}
where we in the last equality have defined $s_k \triangleq b(\theta)-A x_k$ which is the primal slack of the current iterate, and $\sigma_k\triangleq A p_k$ which is how much the step $p_k$ affects the primal feasibility. With this notation, $x_k^*$ being primal feasible is equivalent to $s^*_k \geq 0$.  

If $x^*_k$ is primal infeasible, i.e., if $s^*_k \ngeq 0$, there will be at least one hyper-plane corresponding to an inactive constraint that separates $x_k$ and $x^*_k$. The move from $x_k$ to $x^*_k$ cannot, hence, be completed without breaking feasibility. Instead, a step is taken in the direction of $p_k$ until the first blocking constraint $m\in \bar{\mathcal{W}}_k$ is encountered.  The maximal step length $\alpha_k$ that retains feasibility is explicitly given as 
\begin{equation}
  \label{eq:step-length}
  \alpha_k = \underset{i \in \bar{\mathcal{W}}_k:[s^*_k]_i < 0}{\text{min}} \alpha^i_k, \quad \alpha^i_k \triangleq \frac{[s_k]_i}{[\sigma_k]_i}= \frac{[b]_i + [W]_i \theta - [A]_i x_k}{[A]_i p_k} 
\end{equation}
where $\alpha^j_k$ can be seen as a measure of the distance from the current iterate $x_k$ to the hyper-plane $[A]_j x = [b(\theta)]_j$ in the search direction $p_k$. 

In addition to updating the iterate $x_{k+1} = x_k + \alpha_k p_k$, the working set is updated by adding the first blocking constraint, i.e., the minimizing index of \eqref{eq:step-length}. Concretely, if $m$ is the minimizing index in \eqref{eq:step-length}, the updated working set becomes $ {\mathcal{W}_{k+1} = \mathcal{W}_k \cup \{m\}}$. 

\begin{remark}
  Possible blocking constraints are given by indices in the set $\bar{\mathcal{W}}^{{-}}_k \triangleq {\{i\in \bar{\mathcal{W}}_k : [s^*_k]_i < 0\}}$ since these constraints lead to primal infeasibility when moving from $x_k$ to $x^*_k$ along $p_k$. This is in contrast to the active-set algorithm presented in \cite[Sec. 16-5]{nocedal}, where $\{i\in \bar{\mathcal{W}}_k$: $[\sigma_k]_i > 0\}$ are considered as possible blocking constraints. $[s^*_k]_i <0$ is more restrictive, hence, fewer divisions have to be made in \eqref{eq:step-length} with the formulation in Algorithm \ref{alg:primal-as}.
\end{remark}

If $x_k^*$ is feasible, i.e. if $s^*_k \geq 0$, global optimality for $x^*_k$ is checked by examining the dual variables $\lambda_k$. $x^*_k$ will ge a global optimum if $\lambda_k$ is dual feasible, i.e., if  ${[\lambda_k]_i \geq -\epsilon_d,\:  \forall i \in \mathcal{W}_k}$, where $\epsilon_d$ is the tolerance for dual feasibility. From the first row in \eqref{eq:kkt}, $\lambda_k$ can be obtained by solving 
\begin{equation}
  \label{eq:dual_var}
  A_k^T \lambda_k = -(H x^*_k + f(\theta))
\end{equation}
If the dual iterate is not dual feasible, a constraint corresponding to the most negative dual variable $[\lambda_k]_l$ is removed from the working set, resulting in $\mathcal{W}_{k+1} = \mathcal{W}_k \setminus \{l\}$. 

After the working set has been updated, a new search direction is computed by solving \eqref{eq:kkt-step} with the new working set and the algorithm reiterates the steps described above until global optimality is ensured. 

\begin{remark}
  A straightforward way for terminating the algorithm earlier is to increase $\epsilon_d$, which is further explored in \cite{ifac2020}. 
\end{remark}
\section{Properties of primal active-set algorithms}
\label{sec:prop}
The main operations of Algorithm \ref{alg:primal-as} are removing and adding constraints to the working set. We now consider properties of subsequent search directions and iterates after constraints are added to $\mathcal{W}$, discussed in \ref{sssec:add-proper}, and removed from $\mathcal{W}$, discussed in  \ref{sssec:remove-proper}. 
These insights will later be used in Section \ref{sec:cert-primal-as} to certify the complexity of Algorithm \ref{alg:primal-as}.
\subsection{Addition of a constraint to $\mathcal{W}$}
\label{sssec:add-proper}
When a constraint is added to $\mathcal{W}$ there will be a relationship between the subsequent and previous search direction in terms of $H^*$, as is shown in the following lemma 
\begin{lemma}
  \label{prop:add}
  If a constraint is added to $\mathcal{W}$ in iteration $k$,
  ${p_{k+1} = (1-\alpha_k) H^*_k H p_k}$ 
\end{lemma}
\begin{IEEEproof}
  From the KKT-conditions we have 
  \begin{subequations}
  \begin{align}
	\label{eq:proof-add-stat} H x^*_{k+1} + A^T_{k+1} \lambda_{k+1} &= -f(\theta) \\ 
	\label{eq:proof-add-primal} A x^*_{k+1} &= b_k(\theta)
  \end{align}
  \end{subequations}
  Subtracting $H x_{k+1}$ from \eqref{eq:proof-add-stat} gives 
\begin{equation}
  \begin{aligned}
	\label{eq:proof-mod-stat}
	H p_{k+1} + A_{k+1}^T \lambda_{k+1} &= -f(\theta)-H x_{k+1}   \\ 
										&=-f(\theta)-H x^*_{k} +(1 \text{-} \alpha_{k}) H p_{k}\\
										&= A^T_{k} \lambda_{k} + (1-\alpha_{k}) H p_{k} 
  \end{aligned}
\end{equation}
where $x_{k+1} = x^*_k -(1-\alpha_k)p_k$ has been used in the second equality and $H x^*_{k} + A^T_k \lambda_k = - f(\theta)$ has been used in the third equality. Furthermore, subtracting $A x_{k+1}$ from \eqref{eq:proof-add-primal} gives 
\begin{equation}
  \label{eq:proof-mod-primal}
  A_{k+1} p_{k+1} = b_{k+1} -A_{k+1} x_{k+1} = 0
\end{equation} 
where the last equality follows since $x_{k+1} \in P_{k+1}$. Combining \eqref{eq:proof-mod-stat} with \eqref{eq:proof-mod-primal} gives the KKT-system 
\begin{equation}
  \label{eq:kkt-add}
\begin{pmatrix}
  H & A^T_{k+1}\\
  A_{k+1} & 0 
\end{pmatrix}
\begin{pmatrix}
  p_{k+1} \\ 
  \tilde \lambda
\end{pmatrix} = 
\begin{pmatrix}
  (1-\alpha_k) H p_k \\
 0
\end{pmatrix}
\end{equation}
with $\tilde \lambda = \lambda_{k+1} -  \begin{pmatrix}
  \lambda_k^T & 
0
\end{pmatrix}^T$.
Equation system \eqref{eq:kkt-add} is in the form of the KKT-system in \eqref{eq:kkt}. Hence, by setting ${f(\theta) = \text{-}(1-\alpha_k)H p_k}$, ${b_k(\theta) = 0,}$ and $ x^*_k = p_{k+1}$ in \eqref{eq:kkt-explicit} gives ${p_{k+1} = (1-\alpha_k) \affop_k H  p_k}$.
\end{IEEEproof}

The projective property of $H^*$ from Lemma \ref{lem:proj-operator} can be used together with Lemma \ref{prop:add} to establish a relationship between search directions when constraints are added in consecutive iterations 
\begin{corollary}
  \label{prop:add-col}
  If constraints are added to $\mathcal{W}$ from iteration $k$ until iteration $k+N$, $p_{k+N} = (1-\tau)\affop_{k+N} H p_k$ for some $\tau \in [0,1)$. 
\end{corollary}
\begin{IEEEproof}
  By recursively applying Lemma \ref{prop:add} we get 
  \begin{equation*}
	\begin{aligned}
	p_{k+N} &= \prod_{i=k+N-1}^{k}\Big((1-\alpha_{i})\affop_{i} H \Big) p_k = (1-\tau) \affop_{k+N} H p_k     
	\end{aligned}
  \end{equation*}
  with $(1-\tau) \triangleq \prod_{i=k}^{k+N-1}(1-\alpha_{i})$. The last equality follows from Lemma \ref{lem:proj-operator}, i.e., $\affop_{k+1} H \affop_{k} = H^*_{k+1}$ if $P_{k+1} \subseteq P_k$. Finally, $\tau \in [0,1)$ follows from $\alpha_i \in [0,1), \forall i\in\{k,k+1,\dots,k+N\}$ since constraints were added from iteration $k$ until iteration $k+N$.
\end{IEEEproof}

Corollary \ref{prop:add-col} can be used to get an explicit expression of $x_k$ in terms of $x_0$ and $p_0$ if only additions of constraints have been made since the start of Algorithm \ref{alg:primal-as} up until iteration $k$. 
\begin{corollary}
  \label{prop:add-col2}
  If constraints are added to $\mathcal{W}$ from iteration $0$ until iteration $k$, $x_{k} = \affop_k H(x_0+\tau p_0) +\affconst_k b(\theta)$ 
\end{corollary}
\begin{IEEEproof}
  Using Corollary \ref{prop:add-col} gives 
  \begin{equation*}
  \begin{aligned}
	\tau\affop_k H  p_0 &=  \affop_k H p_0 - p_k = \affop_k H (x^*_0 - x_0) - (x^*_k - x_k) \\    
						&= \affop_k H ( -\affop_0 f(\theta) + \affconst_0 b(\theta)) - \affop_k H x_0 \\& \quad \:\: + \affop_k f(\theta) - \affconst_k b(\theta) + x_k \\ 
						&= -\affop_k H x_0 - \affconst_k b(\theta) + x_k \\
						&\iff x_k = \affop_k H (x_0 + \tau p_0)+ \affconst_k b(\theta)
  \end{aligned}
  \end{equation*}
  where \eqref{eq:kkt-explicit}  has been used in the third equality and the fourth equality follows from Lemma \ref{lem:proj-operator} and \eqref{eq:Hstar-def}.
% H* H T = 0
\end{IEEEproof}

\subsection{Removal of a constraint from $ \mathcal{W} $}
\label{sssec:remove-proper}
When a constraint is removed there will be a relationship between the subsequent search direction and the normal of the removed half-plane, as described by the following lemma 
\begin{lemma}
  \label{prop:remove}
  If constraint $l$ is removed from $\mathcal{W}_k$ in iteration $k$, $p_{k+1} = -[\lambda_k]_l \affop_{k+1} [A]^T_l$
\end{lemma}
\begin{IEEEproof}
  A constraint is removed from $ \mathcal{W}_k$ when a constrained stationary point has been reached. Thus, ${x_{k+1} = x_k^*}$ and the search direction is given by 
\begin{equation}
  p_{k+1} = x^*_{k+1} - x_{k+1} = x^*_{k+1} - x^*_k 
\end{equation}
Since $x^*_{k+1}$ and $x^*_{k}$ are optimal, the following equations hold from the KKT-conditions 
\begin{subequations}
  \begin{align}
  \label{eq:m_kkt_stat}
	H x_k^* + A^T_{k+1} [\lambda_k]_{\bar{l}} + [A]^T_l [\lambda_k]_l &= -f \\ 
  \label{eq:m_kkt_primal}
	A_{k+1} x^*_k &= b_{k+1} \\ 
  \label{eq:m1_kkt_stat}
	H x_{k+1}^* + A^T_{k+1} ([\lambda_k]_{\bar{l}}+\Delta_{\lambda}) &= -f \\  
  \label{eq:m1_kkt_primal}
	A_{k+1} x^*_{k+1} &= b_{k+1} 
  \end{align}
\end{subequations}
where $[.]_{\bar{l}}$ denotes all rows except the $l$:th row. By subtracting \eqref{eq:m_kkt_stat} from \eqref{eq:m1_kkt_stat} and \eqref{eq:m_kkt_primal} from \eqref{eq:m1_kkt_primal} the following KKT-system is obtained 
\begin{equation}
\begin{pmatrix}
  H & A^T_{k+1}\\
  A_{k+1} & 0 
\end{pmatrix}
\begin{pmatrix}
  p_{k+1} \\ 
  \Delta_{\lambda}
\end{pmatrix}
= 
\begin{pmatrix}
  -[A]_l^T [\lambda_k]_l\\
 0 
\end{pmatrix}
\end{equation}
which is in the form of \eqref{eq:kkt} by setting $f(\theta) = [A]^T_l[\lambda_k]_l$, $b(\theta) = 0$ and $x_k^* = p_{k+1}$. Inserting this in \eqref{eq:kkt-explicit} gives 
\begin{equation}
  p_{k+1} =   -\affop_{k+1} ([A]^T_l [\lambda_k]_l) = -[\lambda_k]_l \affop_{k+1} [A]^T_l
\end{equation}
%or in terms of the reduced hessian by 
%\begin{equation}
%  p_{k+1} =  \lambda_k^l Z_{k+1}\inv{(Z_{k+1}^T H Z_{k+1})}Z_{k+1}^T a_l   
%\end{equation}
%with $\gamma \triangleq [\lambda_k]_l$. Since the dual variable of a removed constraint is negative, $\gamma  < 0$ 
which is the stated relation.
\end{IEEEproof}

Lemma \ref{prop:remove} together with Corollary \ref{prop:add-col} gives the following fundamental property of the search directions computed by Algorithm \ref{alg:primal-as} 
\begin{corollary}
  \label{prop:remove-col}
  At iteration $k+N$, let $l$ be the index of the latest removed constraint from $\mathcal{W}$, removed in iteration $k$. Then $p_{k+N} = -(1-\tau) [\lambda_k]_l \affop_{k+N} [A]^T_l$ for some $\tau \in [0,1]$.
\end{corollary}
\begin{IEEEproof}
  The corollary follows from combining Corollary \ref{prop:add-col}, Lemma \ref{prop:remove} and Lemma \ref{lem:proj-operator}. 
\end{IEEEproof}

In conclusion, the search directions will be \textit{completely} determined by $H^*$ acting on the normal of the latest constraint removed from $\mathcal{W}$. Also note that a consequence of this is that the parameter $\theta$ does not affect the direction of the step, only the scaling. This property will be important in the certification of Algorithm \ref{alg:primal-as}, presented in the next section.

\section{CERTIFICATION OF ACTIVE-SET METHOD}
\label{sec:cert-primal-as}
This section describes a method to exactly identify which sequence of working-set changes different parameters will give rise to when Algorithm \ref{alg:primal-as} is applied to \eqref{eq:mpQP}. For the time being we assume, for clarity, that the reduced Hessian in nonsingular, i.e., that the KKT-system \eqref{eq:kkt} has a unique solution. In Section \ref{sec:singular} we amend the method for the singular case.  The method is an extension of \cite{cdc2019} and similar to the ones presented in \cite{10.1109TAC.2017.2696742}, \cite{cimini2019complexity}, and \cite{10.1109TAC.2011.2108450}, in the sense that the parameter space is iteratively partitioned depending on how the working set changes in each iteration.

There are two sources leading to a change in the working set: either a constraint is added or removed. A removal only happens after a constrained stationary point has been reached. Moreover, if this point is a global optimum, i.e., if all the dual variables are nonnegative, Algorithm \ref{alg:primal-as} terminates with the global solution. In contrast, a constraint will be added to $\mathcal{W}$ if there is a blocking constraint between the current iterate and constrained point. Thus, Algorithm \ref{alg:primal-as} can be split into two modes 
\begin{enumerate}[a)]
  \item Checking for global optimality and removing constraints, performed at lines  6-9.
  \item Checking for local optimality and adding constraints, performed at lines 3-5 and 9-11.  
\end{enumerate}
The algorithm goes from mode a) $\rightarrow$ b) when a constraint is removed, whereas it goes from mode b) $\rightarrow$ a) when a constrained stationary point is primal feasible. 
\usetikzlibrary{shapes,arrows}

% Define block styles
\tikzstyle{decision} = [diamond, draw, fill=white!20, 
text width=4.5em, text badly centered, node distance=2cm, inner sep=0pt]
\tikzstyle{block} = [rectangle, draw, fill=white!20, 
text width=10em, text centered, rounded corners, minimum height=4em]

 Single column
\tikzstyle{decision} = [diamond, draw, fill=white!20, 
text width=4.5em, text badly centered, node distance=2.5cm, inner sep=0pt]
\tikzstyle{block} = [rectangle, draw, fill=white!20, 
text width=10em, text centered, rounded corners, minimum height=4em,node distance=2.5cm]

\tikzstyle{line} = [draw, -latex']
\begin{figure}[htpb]
\begin{center}
\begin{tikzpicture}[scale=0.8, transform shape,node distance=2.25cm, auto]
  % Place nodes
  \node [block,line width=0mm,draw=gray!60] (init) {Starting iterate $x_0$ and working set $\mathcal{W}_0$};
  \node [block,fill=red!20,line width=0.35mm, below of=init, node distance=1.7cm] (add1) {Compute $x^*_k$};
  \node [decision, below of=add1, fill=red!20,line width=0.35mm] (CSP) {$s^*_k \geq 0$ };
  \node [block,fill=green!20,line width=0.35mm, below of=CSP] (lambda) {Compute $\lambda_k$};
  \node [decision,fill=green!20,line width=0.35mm, below of=lambda] (dual) {$\lambda_k \geq \text{-}\epsilon_d$};
  \node [block, line width=0.35mm, fill=green!20, right of=dual,node distance=4.1cm] (remove) {Remove constraint from $\mathcal{W}$ and update $x$};
  \node [block, below of=dual,line width=0mm,draw=gray!60] (end) {$x^*$ found};
  \node [block,fill=red!20,line width=0.35mm, right of=CSP,node distance=4.1cm] (add2) {Add constraint to $\mathcal{W}$ and update $x$};
  \node [right of=remove] (rm_help)  {};
  % Draw edges
  \path [line] (init) -- (add1);
  \path [line] (add1) -- (CSP);
  \path [line] (CSP) -- (lambda);
  \path [line] (lambda) -- (dual);
  \path [line] (dual) --node[right]{Yes} (end);
  \path [line] (dual) --node[above]{No} (remove);
  \path [line] (CSP) --node[right]{Yes} (lambda);
  \path [line] (CSP) --node[above]{No} (add2);
  \path [line] (add2) |- (add1); 
  %\path [line] (remove) -- ++(2.25,0) |-(add1); 
  %Single column
  \path [line] (remove) -- ++(2.75,0) |-(add1); 
  % Legend
  \node[rectangle,fill=green!20,minimum width=30mm] (aff) at (4,0.25){Mode a)}; 
  \node[rectangle,fill=red!20,minimum width=30mm] (quad) at (4,-0.3){Mode b)}; 
\end{tikzpicture}
\end{center}
\caption{Flowchart characterizing Algorithm \ref{alg:primal-as}}
\label{fig:main-alg}
\end{figure}
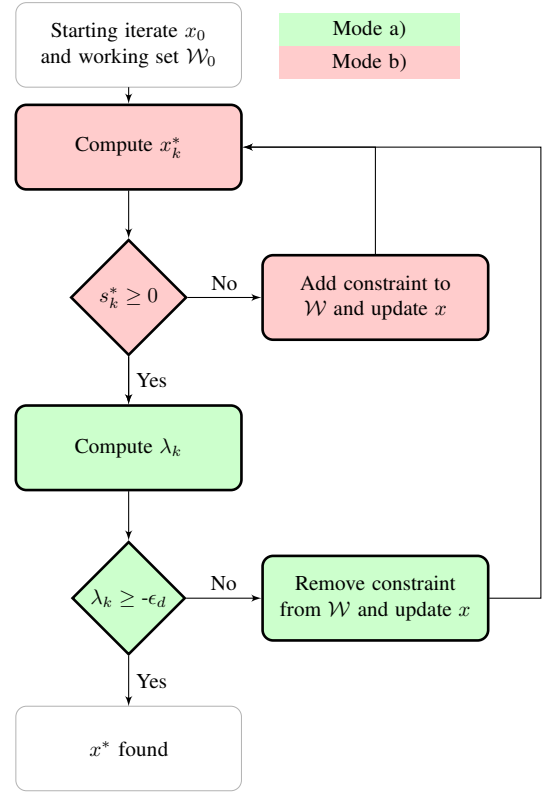
This characterization of Algorithm \ref{alg:primal-as}, illustrated in Figure \ref{fig:main-alg}, is used to create a partition of $\Theta_0$ reflecting which sequence of working-set changes different parameters generate when Algorithm \ref{alg:primal-as} is applied to the mpQP in \eqref{eq:mpQP}. 

Parameter regions are partitioned in the following way:   
If a region $\Theta$ is in mode a) it will be partitioned into the following parameter regions 
\begin{itemize}
  \item $\Theta^*$ -  Global optimality obtained.
  \item $\Theta^{\text{-}j}$ - $j$ removed from $\mathcal{W}$.
\end{itemize}
Likewise, in mode b) a region $\Theta$ will be partitioned into the following parameter regions 
\begin{itemize}
  \item $\Theta^{\text{CSP}}$ - A primal feasible CSP reached.
  \item $\Theta^{\text{+}j}$ - $j$ added to $\mathcal{W}$. 
\end{itemize}
$\Theta_0$ will iteratively be partioned into these subsets, corresponding to executing iteartions of Algorithm \ref{alg:primal-as} parametrically, until all parameters have reached global optimality. In the final partition, parameters in the same region signify that they produce the same sequence of working-set changes to reach optimality. The method is summarized in Algorithm \ref{alg:main}.

Each region of the partition is represented by a tuple $(\Theta,\mathcal{W},F,G,s, k, \hat{n})$ containing the following data  
\begin{itemize}
  \item $\Theta \subseteq \Theta_0 \subseteq \mathbb{R}^{p}$ -  The subset of the parameter space that defines the region.
  \item $ \mathcal{W}$ - The working set in the region, for the current iteration. 
  \item $F\in \mathbb{R}^{n\times p} $ and $G \in \mathbb{R}^{n\times1}$ - Matrices that define the affine mapping $x_k = F \theta + G$ for $\theta \in \Theta$. 
  \item $s$ - A status flag that marks if the region has reached a CSP 1, globally optimality 2, or neither 0. 
  \item $k$ - Number of iterations performed by Algorithm \ref{alg:primal-as} to reach the current state.  
  \item $\hat{p}$ - The normal of the latest constraining half-plane that has been removed from the working set.
\end{itemize}

$S$ is a stack containing tuples corresponding to regions of $\Theta_0$ that are yet to reach global optimality.
\begin{remark}
  Algorithm \ref{alg:main} is well suited for parallelization by distributing the stack $S$ over multiple processors.   
\end{remark}
\begin{algorithm}[H]
  \caption{Partition $\Theta_0$ based on working-set changes}\label{alg:main}
  \begin{algorithmic}[1]
	\Require $\Theta_0, \mathcal{W}_0, F_0, G_0,$ mpQP 
	\Ensure FinalPartition
	\State Push $(\Theta_{0}, \mathcal{W}_0,F_0,G_0,0,0,\text{NaN})$ to $S$
  \While {$S$ is not empty}
  \State Pop $p_c$ from $S$
  %\If {\textsc{isCSP}$(p_c, \text{mpQP})$}
  \If {$p_c$ has reached a CSP}
  \State Partition = \textsc{ModeA}($p_c$,mpQP)
  \Else
  %\If {Reduced Hessian for $p_c$ is nonsingular} 
  \State Partition = \textsc{ModeB}($p_c$, mpQP)
  %\Else 
  %\State Partition = \textsc{SingularModeB}($p_c$, mpQP)
  %\EndIf
  \EndIf
  \For {$p$ in Partition}
  \If {$p$ is global optimum}
  \State Append $p$ to FinalPartition
  \Else 
  \State Push $p$ to $S$
  \EndIf
  \EndFor
  \EndWhile
  \State \Return FinalPartition
  \end{algorithmic}
\end{algorithm}

In Algorithm \ref{alg:main}, the procedure \textsc{ModeA} partitions the parameter space depending on what happens in mode a), i.e., whether global optimality is reached or if a constraint is removed. The procedure is described in detail in Section \ref{ssec:remove} and is summarized in Algorithm \ref{alg:remove} in the end of that section. 
Likewise, the procedure \textsc{ModeB} partitions the parameter space depending on what happens in mode b), i.e., whether a CSP is reached or if a constraint is added. The procedure is described in detail in Section \ref{ssec:add} and is summarized in Algorithm \ref{alg:add} in the end of that section.

\subsection{Removing constraints and checking for global optimality}
\label{ssec:remove}
How the parameter space is partitioned in mode a) will now be described in detail. 
At iteration $k$, the variable that decides whether global optimality has been reached or if a constraint has to be removed is $\lambda_k(\theta)$. Recall from Algorithm \ref{alg:primal-as} that a global optimum has been found at iteration $k$ if all $\lambda_k(\theta)$ are nonnegative, within a given tolerence $\epsilon_d$. Otherwise, a constraint $l$ corresponding to a negative dual-variable is removed from the working set. From Algorithm \ref{alg:primal-as} line \ref{line:const-rm}, $l$ is chosen as the most negative component of $\lambda_k(\theta)$, i.e.,    
\begin{equation}
  l = \underset{i \in \mathcal{W}_k}{\text{argmin}} [\lambda_k(\theta)]_i.
\end{equation}
Hence, the set $\Theta_k^{\text{-}j}$ of all parameters in iteration $k$ resulting in constraint $j\in \mathcal{W}_k $ being removed from the working set is given by
\begin{equation}
  \begin{aligned}
  \label{eq:th-remove}
	\Theta_k^{\text{-}j} = \{ \theta \in \Theta_k | &[\lambda_k(\theta)]_j < \text{-} \epsilon_d\\ &[\lambda_k(\theta)]_j < [\lambda_k(\theta)]_i,\:\: \forall i \in \mathcal{W}_k \setminus \{j\}\} 
  \end{aligned}
\end{equation}
i.e., $\theta$ for which the dual variable corresponding to constraint $j$ is negative and more negative than any other dual variable. 

Likewise, the set $\Theta_k^*$ of all parameters in iteration $k$ resulting in a global optimum is given by 
\begin{equation}
  \label{eq:th-global}
  \Theta_k^* = \{\theta \in \Theta_k | [\lambda_k(\theta)]_i \geq \text{-}\epsilon_d,\:\: \forall i \in \mathcal{W}_k \} 
\end{equation}
i.e., $\theta$ for which all dual variables are nonnegative.

To summarize, a region $\Theta_k$ will be partitioned into $\Theta^*_k$ and $\Theta^{\text{-}i}_k, \forall i\in {\mathcal{W}}_k$ in mode a), as illustrated in Figure \ref{fig:modeb}. 
\usetikzlibrary{shapes,arrows,calc}
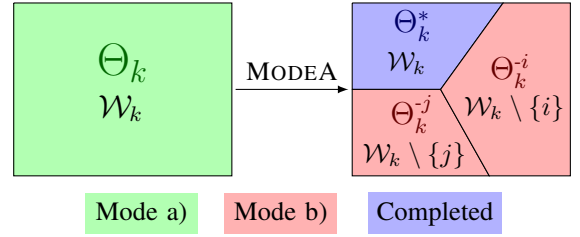
\begin{figure}[htpb]
  \begin{center}
	\begin{tikzpicture}[scale=1, transform shape]
	  \coordinate(s1) at (-5.55,1.15);
	  \coordinate(s2) at (-2.65,1.15);
	  \coordinate(s3) at (-2.65,-1.15);
	  \coordinate(s4) at (-5.55,-1.15);

	  \coordinate(v1) at (-1.05,1.15);
	  \coordinate(v2) at (1.85,1.15);
	  \coordinate(v3) at (1.85,-1.15);
	  \coordinate(v4) at (-1.05,-1.15);
	  \coordinate(u1) at (1.75,2.25);
	  \coordinate(u2) at (-1.5,-2.25);
	  \coordinate(w1) at (-2,0);
	  \coordinate(w2) at (3.5,0);
	  \coordinate(q1) at (-1,2);
	  \coordinate(q2) at (1.25,-2);
	  \coordinate (key) at (intersection of u1--u2 and w1--w2) {};
	  \coordinate (optv1) at (intersection of u1--u2 and v1--v2){};
	  \coordinate (optv2) at (intersection of w1--w2 and v4--v1){};
	  \coordinate (add1) at (intersection of key--q2 and v4--v3){};

	  \draw[fill,blue!30!white,draw=black] (v1)--(optv1)--(key)--(optv2)--(v1);
	  \draw[fill,red!30!white,draw=black] (optv2)--(key)--(add1)--(v4)--(optv2);
	  \draw[fill,red!30!white,draw=black] (optv1)--(v2)--(v3)--(add1)--(key)--(optv1);
	  \draw[fill,green!30!white,draw=black] (s1)--(s2)--(s3)--(s4)--(s1);
	  \node[red!40!black] (jremove) at (-0.255,-0.33) {\large$\Theta_k^{\text{-}j}$};
	  \node[] (jremove) at (-0.225,-0.875) {$\mathcal{W}_k\setminus\{j\}$};
	  %\node[red!40!black] (jremove) at (0.05,-0.85) {$\mathcal{W}_{k\text{+}1} = \mathcal{W}_k \cup \{ j\}$};
	  \node[red!40!black] (iremove) at (1.05,0.25) {\large$\Theta_k^{\text{-}i}$};
	  \node[] (iremove) at (1.1,-0.25) {$\mathcal{W}_k\setminus\{i\}$};

	  \node[blue!40!black] (star) at (-0.25,0.85) {\large$\Theta_k^*$};
	  \node[] (star) at (-0.325,0.35) {$\mathcal{W}_k$};

	  \node[green!40!black] (start) at (-4.1,0.35) {\LARGE$\Theta_k$};
	  \node[] (start) at (-4.125,-0.25) {\large$\mathcal{W}_k$};

	  \draw[-latex] (-2.6,0) -- node[above]{\textsc{ModeA}} (-1.1,0);
	  \node[rectangle,fill=green!30!white,minimum width=15mm] (aff) at (-3.85,-1.65){Mode a)}; 
	  \node[rectangle,fill=red!30!white,minimum width=15mm] (aff) at (-2,-1.65){Mode b)}; 
	  \node[rectangle,fill=blue!30!white,minimum width=15mm] (aff) at (.05,-1.65){Completed}; 
  \end{tikzpicture}
\end{center}
\caption{Partitioning of a region $\Theta_k$ performed in mode a).}
\label{fig:modea}
\end{figure}
To get more explicit expressions of these sets, recall from \eqref{eq:kkt_para_sol} that $\lambda_k(\theta)$ is affine in $\theta$, i.e., ${\lambda_k(\theta) = F_k^{\lambda}  \theta + G_k^{\lambda}}$.
Using this, the regions $\Theta_k^{\text{-}j}$ in \eqref{eq:th-remove} can be equivalently expressed as all $\theta \in \Theta_k$ such that  
\begin{subequations}
  \label{eq:th-remove}
  \begin{align}
	[F^{\lambda}_k]_j\theta + [G^{\lambda}_k]_j &< \text{-} \epsilon_d\\ 
	([F^{\lambda}_k]_j\text{-}[F^{\lambda}_k]_i)\theta &< ([G^{\lambda}_k]_i\text{-}[G^{\lambda}_k]_j), \:\: \forall i \in \mathcal{W}_k \setminus \{j\} 
  \end{align}
\end{subequations}
Likewise, the region $\Theta_k^*$ defined in \eqref{eq:th-global} can be equivalently expressed as 
\begin{equation}
  \label{eq:th_global_opt}
  \Theta_k^* = \{\theta \in \Theta_k |F_k^{\lambda}\theta + G_k^{\lambda} \geq \text{-} \epsilon_d \} 
\end{equation}
How regions of the parameter space are partitioned in mode a) is summarized in Algorithm \ref{alg:remove}.
\begin{remark}
  \label{rem:case1part}
  Importantly, all partitioning in \eqref{eq:th-remove} and \eqref{eq:th_global_opt} are made by linear inequalities. 
\end{remark}
\begin{algorithm}[H]
  \caption{Partition $\Theta$ based on if global optimality is reached or if a constraint is removed from $\mathcal{W}$ }\label{alg:remove}
  \begin{algorithmic}[1]
	%\Require $(\Theta,\mathcal{A},F,G,s, k, \hat{p})$ , mpQP
	%\Ensure $P$
	\State \textsc{ModeA}{($(\Theta,\mathcal{W},F,G,s, k, \hat{n})$ , mpQP)}
	%\State $\verb!++!k$
	\State Calculate $F^{\lambda}$ and $G^{\lambda}$ according to \eqref{eq:F-lam}
	\For {all $i$ in $ \mathcal{W}$}
	\State Calculate $\Theta^{\text{-}i}$ according to \eqref{eq:th-remove} 
	\If {$\Theta^{\text{-}i} \neq \emptyset $ }  \label{line:alg4-empty1}
	\State Append$\:\:(\Theta^{\text{-}i}, \mathcal{W}{\setminus}\{i\},F,G,0,k\text{+}1, [A]^T_i)$ to $P$ 
	\EndIf
	\EndFor
	\State Calculate $\Theta^*$ according to \eqref{eq:th_global_opt} 
	\If {$\Theta^* \neq \emptyset $ } \label{line:alg4-empty2}
	\State Append $(\Theta^*, \mathcal{W},F,G,2,k, \hat{n})$ to $P$ 
	\EndIf
	\State \Return $P$
  \end{algorithmic}
\end{algorithm}

\subsection{Adding constraints and checking for local optimality}
\label{ssec:add}

We now turn our attention to how the parameter space is partitioned in mode b). 
If $j$ is the minimizing index of the minimization in \eqref{eq:step-length}, it will be added to $\mathcal{W}_{k+1}$ and $\alpha_k = \alpha^j_k$.
The set $\Theta_k^{\text{+}j}$ of all parameters in iteration $k$ leading to constraint $j$ being added to $ \mathcal{W}_{k+1}$ is, hence, given by   
\begin{equation}
  \begin{aligned}
  \label{eq:th-add}
\Theta_k^{\text{+}j} \triangleq \{\theta \in \Theta_k| [s_k^*]_j < 0,\: \alpha_k^j(\theta) < \alpha_k^i(\theta), \forall i \in \bar{\mathcal{W}}_k^{-}\setminus\{j\}\} 
  \end{aligned}
\end{equation}
where $j$ being a blocking constraint is ensured by $[s^*_k]_j < 0$, while $\alpha_k^j(\theta) < \alpha^i_k(\theta),\:\: \forall i \in \bar{\mathcal{W}}^-_k \setminus \{j\}$ ensures that it is the \textit{first} encountered blocking constraint.   

Furthermore, the constrained stationary point is primal feasible if $[s^*_k]_i \geq 0$ $\forall i \in \bar{\mathcal{W}}_k$. The set $\Theta_k^{\text{CSP}}$ of all parameters in iteration $k$ leading to a constrained stationary point being reached is, hence, given by
\begin{equation}
  \label{eq:th-local}
  \begin{aligned}
  \Theta_k^{\text{CSP}} \triangleq \{\theta \in \Theta_k |& [s^*_k]_i \geq 0,  &\forall i \in \bar{\mathcal{W}}_k\}  \end{aligned}
\end{equation}

To summarize, a region $\Theta_k$ will be partitioned into $\Theta^{\text{CSP}}_k$ and $\Theta^{\text{+}i}_k, \forall i\in \bar{\mathcal{W}}_k$, in mode b), as illustrated in Figure \ref{fig:modea}. 
\usetikzlibrary{shapes,arrows,calc}
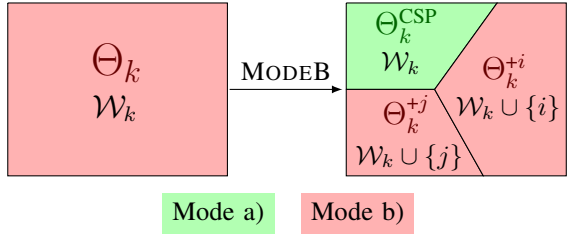
\begin{figure}[htpb]
  \begin{center}
	  \begin{tikzpicture}[scale=1, transform shape]
		\coordinate(s1) at (-5.55,1.15);
		\coordinate(s2) at (-2.63,1.15);
		\coordinate(s3) at (-2.63,-1.15);
		\coordinate(s4) at (-5.55,-1.15);
		
		\coordinate(v1) at (-1.05,1.15);
		\coordinate(v2) at (1.875,1.15);
		\coordinate(v3) at (1.875,-1.15);
		\coordinate(v4) at (-1.05,-1.15);
		\coordinate(u1) at (1.75,2.25);
		\coordinate(u2) at (-1.5,-2.25);
		\coordinate(w1) at (-2,0);
		\coordinate(w2) at (3.5,0);
		\coordinate(q1) at (-1,2);
		\coordinate(q2) at (1.25,-2);
		\coordinate (key) at (intersection of u1--u2 and w1--w2) {};
		\coordinate (optv1) at (intersection of u1--u2 and v1--v2){};
		\coordinate (optv2) at (intersection of w1--w2 and v4--v1){};
		\coordinate (add1) at (intersection of key--q2 and v4--v3){};
		
		\draw[fill,green!30!white,draw=black] (v1)--(optv1)--(key)--(optv2)--(v1);
		\draw[fill,red!30!white,draw=black] (optv2)--(key)--(add1)--(v4)--(optv2);
		\draw[fill,red!30!white,draw=black] (optv1)--(v2)--(v3)--(add1)--(key)--(optv1);
		\draw[fill,red!30!white,draw=black] (s1)--(s2)--(s3)--(s4)--(s1);
		\node[red!40!black] (jremove) at (-0.25,-0.325) {\large$\Theta_k^{\text{+}j}$};
		\node[] (jremove) at (-0.2125,-0.9) {$\mathcal{W}_k\cup\{j\}$};
		%\node[red!40!black] (jremove) at (0.05,-0.85) {$\mathcal{W}_{k\text{+}1} = \mathcal{W}_k \cup \{ j\}$};
		\node[red!40!black] (iremove) at (1.05,0.25) {\large$\Theta_k^{\text{+}i}$};
		\node[] (iremove) at (1.1,-0.25) {$\mathcal{W}_k\cup\{i\}$};
		
		\node[green!40!black] (star) at (-0.25,0.85) {\large$\Theta_k^{\text{CSP}}$};
		\node[] (star) at (-0.325,0.35) {$\mathcal{W}_k$};
		
		\node[red!40!black] (start) at (-4.1,0.35) {\LARGE$\Theta_k$};
		\node[] (start) at (-4.125,-0.25) {\large$\mathcal{W}_k$};

		\draw[-latex] (-2.6,0) -- node[above]{\textsc{ModeB}} (-1.1,0);
	  \node[rectangle,fill=green!30!white,minimum width=15mm] (aff) at (-2.75,-1.65){Mode a)}; 
	  \node[rectangle,fill=red!30!white,minimum width=15mm] (aff) at (-0.9,-1.65){Mode b)}; 
  \end{tikzpicture}
\end{center}
\caption{Partitioning of a region $\Theta_k$ performed in mode b).}
\label{fig:modeb}
\end{figure} 

In the rest of this section, we derive explicit expressions for $\Theta^{\text{CSP}}_k$ and $\Theta^{\text{+}j}_k$. First, we formulate an explicit expression for $\Theta^{\text{CSP}}_k$, which is straightforward since $s^*_k$ is affine in $\theta$, i.e., 
\begin{equation}
  \label{eq:opt-slack}
  s^*_k = F^*_{s_k} \theta + G^*_{s_k}\\
\end{equation}
\begin{equation}
  F^*_{s_k} \triangleq W-AF^*_k, \quad G^*_{s_k} \triangleq b-AG^*_k
\end{equation}
An explicit expression for $\Theta^{\text{CSP}}_k$ is, hence, all ${\theta \in \Theta_k}$ such that  
\begin{equation}
  \label{eq:opt-slack-part}
  [F^*_{s_k}]_i \theta + [G^*_{s_k}]_i \geq 0,\quad \forall i \in \bar{\mathcal{W}}_k
\end{equation}

Next, we formulate an explicit expression for $\Theta^{\text{+}j}_k$. This entails some technicalities which stem from the behaviour of Algorithm \ref{alg:primal-as} being different depending on if a constraint has been removed from $\mathcal{W}$ or not, as was discussed in the end of Section \ref{sec:prop}. Essentially, the analysis becomes simpler after a constraint has been removed from $\mathcal{W}$ since additional structure is introduce to the search direction $p_k$ and, as will be shown, to the iterates $x_k$. Therefore, two different cases are considered when describing $\Theta_k^{\text{+}j}$ explicitly in terms of $\theta$: Case 1 considers the case when a constraint has been removed from $\mathcal{W}$ in an earlier iteration, whereas Case 2 considers the case when no constraint has been removed since the start of Algorithm \ref{alg:primal-as}. 

From \eqref{eq:th-add}, the quantities that define $\Theta^{\text{+}j}_k$ are $s^*_k(\theta)$ and $\alpha^i_k(\theta), i \in \bar{\mathcal{W}}^{-}_k$, where we know from above that $s^*_k(\theta)$ is affine in $\theta$. The main complication for formulatating an explicit expression of $\Theta^{\text{+}j}_k$ is, hence, to establish an explicit expression for $\alpha^i_k(\theta)$, which will be straightforward in Case 1 because of the structure of $p_k$ and $x_k$, and more technical in Case 2. 

\subsubsection{Case 1 - A constraint has been removed from $\mathcal{W}$}
As was mentioned above, the main challenge when expressing an explicit expression for $\Theta_k^{\text{+}j}$ is to express $\alpha^j_k(\theta)$ explicitly which, in turn, requires an explicit expression for the iterate $x_k$ and the search direction $p_k$ since they define $\alpha^j_k$ in \eqref{eq:step-length}. When a constraint has been removed, $p_k$ will, from Corollary \ref{prop:remove-col}, be related to the latest removed constraint $l$, removed in iteration $\tilde k$, by 
\begin{equation}
  \label{eq:nice-step}
  p_k(\theta) = -(1-\tau)[\lambda_{\tilde k}]_l \affop_k [A]_l^T = \gamma(\theta) H^*_k \hat{p}
\end{equation}
with the scaling factor ${\gamma(\theta) \triangleq -(1-\tau)[\lambda_{\tilde k}]_l}$ and the latest removed normal ${\hat{p} \triangleq [A]^T_l}$. Note that $\gamma(\theta)>0$, which follows from $\tau \in[0,1)$ and $[\lambda_{\tilde{k}}]_l<0$ since constraint $l$ was removed in iteration $\tilde{k}$. 

We will now show that the iterates are also endowed with a simple structure after a constraint has been removed from $\mathcal{W}$. Namely,  all subsequent iterates will be affine in $\theta$. 
\begin{theorem}
  \label{prop:affine}
  If a constraint is removed in iteration $\kappa$, ${x_k = F_k \theta + G_k,\:\forall k > \kappa}$ for some $F_k \in \mathbb{R}^{n\times p}, G_k \in \mathbb{R}^{n}$.
\end{theorem}
\begin{IEEEproof}
  Without loss of generality, let $\tilde k \geq \kappa$ be the latest iteration in which a constraint was removed and let $l$ be the corresponding index of the constraint that was removed. Now, assume that $x_k = F_k \theta + G_k$ for $k > \tilde k \geq \kappa$ and first consider the case when there is a blocking constraint. Let $j$ be the corresponding index of the first blocking constraint and let $\hat{p}$ and $\gamma(\theta)$ be defined as above. $x_k = F_k \theta + G_k$ together with the expression of $p_k$ in \eqref{eq:nice-step} inserted into the definition of $\alpha^j_k$ in \eqref{eq:step-length} gives 
  \begin{equation}
  \label{eq:rem-steplength}
  \alpha_k(\theta) = \alpha_k^j (\theta) = \frac{[b]_j + [W]_j \theta - [A]_j (F_k \theta + G_k)}{\gamma(\theta) [A]_j \affop_k \hat{p}} 
\end{equation}
Moreover, recall that the subsequent iterate $x_{k+1}$ is given by
\begin{equation}
  \label{eq:next-iter}
  x_{k+1}(\theta) = x_k(\theta) + \alpha_k(\theta)p_k(\theta)
\end{equation}
By inserting \eqref{eq:rem-steplength} and \eqref{eq:nice-step} in \eqref{eq:next-iter}, after simplifications, one gets $x_{k+1} = F_{k+1} \theta + G_{k+1}$, where $F_{k+1}$ and $G_{k+1}$ are given by 
\begin{equation}
  \label{eq:x_plus}
  \begin{aligned}
	F_{k+1} =& F_{k} + \affop_k \hat p \frac{[W]_j -  [A]_j F_k}{[A]_j \affop_k \hat{p}} \\ G_{k+1} =& G_k + \frac{[b]_j-[A]_j G_k}{[A]_j \affop_k \hat{p}} \affop_k \hat{p}
  \end{aligned}
\end{equation}
If instead there are no blocking constraints $x_{k+1} = x_k^*$, which is affine in $\theta$ by \eqref{eq:kkt_para_sol}, completing the induction step.

Similarly, the base case follows since if a constraint was removed in iteration $\tilde k$, $x_{\tilde k+1} =x^*_{\tilde k}$,  which is affine in $\theta$ by \eqref{eq:kkt_para_sol}. Hence, the theorem follows by induction.
\end{IEEEproof}

With the explicit expression for $x_k$ from Theorem \ref{prop:affine}, and the explicit expression for $p_k$ from \eqref{eq:nice-step},  the step length $\alpha_k^j$, defined in \eqref{eq:step-length}, is given by
\begin{equation}
  \label{eq:alpha_case2_th}
  \alpha_k^j(\theta) = \frac{[F_{s_k}]_j \theta + [G_{s_k}]_j}{\gamma(\theta)[G_{\sigma_k}]_j} 
\end{equation}
with $F_{s_k}, G_{s_k} \text{and } G_{\sigma_k}$ defined as 
\begin{equation*}%
  F_{s_k} \triangleq W - A F_k,\:\: G_{s_k} \triangleq b - A G_k, \:\:G_{\sigma_k} \triangleq A \affop_k \hat{p} 
\end{equation*}
By inserting expression \eqref{eq:alpha_case2_th} for $\alpha^j_k$ and expression \eqref{eq:opt_slack} for $s^*_k$ in \eqref{eq:th-add}, $\Theta^{\text{+}j}_k$ can be explicitly stated as all $\theta \in \Theta_k$ satisfying
\begin{subequations}
  \label{eq:th_add_case2}
  \begin{align}
	&K^{j,i}_k \theta <L^{j,i}_k, \quad \forall i \in \bar{\mathcal{W}}_k\setminus\{ j\}   \\
	&[F^*_{s_k}]_j\theta + [G^*_{s_k}]_j < 0
  \end{align}
\end{subequations}
where $K^{j,i}_k$ and $L^{j,i}_k$ is given by
\begin{subequations}
  \label{eq:KL-def}
  \begin{align}
	K^{j,i}_k &\triangleq [G_{\sigma_k}]_i [F_{s_k}]_j-[G_{\sigma_k}]_j [F_{s_k}]_i\\
	L^{j,i}_k &\triangleq -[G_{\sigma_k}]_i [G_{s_k}]_j+[G_{\sigma_k}]_j [G_{s_k}]_i
  \end{align}
\end{subequations}
\begin{remark}
  \label{rem:case1nice}
  Since all inequalities introduced in Case 1 are affine, see Remark \ref{rem:case1part}, and that Case 2 never occurs again once it has been left - since once a constraint has been removed from $\mathcal{W}$ it is impossible to return to the state of never having removed a constraint - \textit{all further partitioning of the parameter space will exclusively be done by half-planes.}
\end{remark}

\label{sssec:case1}
\subsubsection{Case 2 - No constraint has been removed from $\mathcal{W}$}
\label{sssec:case2}

 When formulating an explicit expressions for $\Theta^{\text{+}j}_k$ when no constraint has been removed from $\mathcal{W}$, we will use the quantity $\tilde\alpha_k^j$ defined as 

\begin{equation}
  \label{eq:alp_tilde}
  \tilde{\alpha}_k^j(\theta) \triangleq \frac{[b]_j + [W]_j \theta- [A]_j (\affop_k H x_0(\theta) + T_k b(\theta))}{[A]_j \affop_k H p_0(\theta)}, 
\end{equation}
instead of $\alpha^j_k$, where $x_0$ is the starting iterate and $p_0(\theta) = x^*_0(\theta) - x_0(\theta)$.
$\tilde{\alpha}_k^j$ can be seen as a measure of the distance between the starting iterate $x_0$ projected onto $P_k$, given by $H^*_k H x_0$, and the half-plane $[A]_j x = [b(\theta)]_j$ along the search direction. Figure \ref{fig:alpha} depicts a simple two-dimensional case to capture the relationship between $\alpha^j_k$ and $\tilde\alpha^j_k$. 
\usetikzlibrary{intersections}
\usetikzlibrary{decorations.pathreplacing}
\usetikzlibrary{calc}
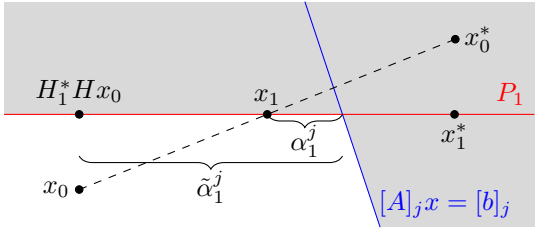
\begin{figure}[htpb]
\begin{center}
\begin{tikzpicture}[scale=1, transform shape]
  \coordinate(P1init) at (-1,1);
  \coordinate(P1end) at (6.05,1);
  \coordinate(Pjinit) at (3,2.5);
  \coordinate(Pjend) at (4,-0.5);
  \coordinate(x0) at (0,0);
  \coordinate(x0star) at (5,2);
  \coordinate (x1) at (intersection of x0--x0star and P1init--P1end) {};
  \coordinate (pix0) at (intersection of x0--++(0,10) and P1init--P1end) {};
  \coordinate (x1star) at (intersection of x0star--{++(0,-500)} and P1init--P1end) {};
  \coordinate (x2) at (intersection of Pjinit--Pjend and P1init--P1end) {};
  \fill[gray!30!white] (P1init) rectangle ($(P1end)+(0,1.5)$);
  \fill[gray!30!white] (Pjend) rectangle (P1end);
  \fill[gray!30!white] (x2)--(Pjend) -- (P1end)--cycle;
  \draw[red] (P1init)-- (P1end)node[above left]{$P_1$}; 
  \draw[blue] (Pjinit) --(Pjend)node[above right,xshift=-4pt]{$[A]_j x = [b]_j$};
%  \draw[dashed] (0,0) node[left]{$x_0$} coordinate (x0)--(5,2) node[right]{$x^*_0$} coordinate (x0star);
  \draw[fill] (x0) node[left]{$x_0$} circle [radius=.05];
  \draw[fill] (x0star) node[right]{$x^*_0$} circle [radius=.05];
  %\draw[dashed,green!40!black] (x0) -- (x0star) node [near end, above] {$L_0$};;
  \draw[dashed] (x0) -- (x0star);
  \draw[fill] (x1) node[above]{$x_1$} circle [radius=.05];
  \draw[fill] (pix0) node[above]{$H^*_1 H x_0$} circle [radius=.05];
  \draw[fill] (x1star) node[below]{$x^*_1$} circle [radius=.05];
  \draw [decorate,decoration={brace,amplitude=5pt,mirror,raise=0ex}]
  (x1) -- (x2) node[midway,yshift=-1em]{$\alpha_1^j$};
  \draw [decorate,decoration={brace,amplitude=5pt,mirror,raise=3.5ex}]
  (pix0) -- (x2) node[midway,yshift=-2.75em]{$\tilde\alpha_1^j$};
  \end{tikzpicture}
\end{center}
\caption{Relationship between $\alpha^j_1$ and $\tilde\alpha^j_1$ for a fixed $\theta$. $\tilde\alpha^j_1$ and $\alpha^j_1$ are fractions of a full step to $x^*_1$, not geometric distances. The white and grey areas mark the feasible set and its complement, respectively.}
\label{fig:alpha}
\end{figure}
The main reason for considering $\tilde{\alpha}^j_k(\theta)$ instead of $\alpha^j_k(\theta)$ is that $\alpha^j_k(\theta)$ dependence on $\theta$ in an intricate way, whereas $\tilde{\alpha}^j_k$ simply is a linear fraction of $\theta$ 
\begin{equation}
  \label{eq:sur-steplength}
  \tilde\alpha_k^j(\theta) = \frac{[\tilde{F}_{s_k}]_j \theta + [\tilde{G}_{s_k}]_j}{[\tilde{F}_{\sigma_k}]_j \theta + [\tilde{G}_{\sigma_k}]_j} 
\end{equation}
with $\tilde{F}_{s_k},\tilde{G}_{s_k},\tilde{F}_{\sigma_k} \text{ and } \tilde{G}_{\sigma_k}$ defined as 
\begin{subequations}
  \begin{align}
	\label{eq:Cj_case1}
	\tilde{F}_{s_k} &\triangleq  W - A (\affop_k H F_0 + \affconst_k W) \\  
	\tilde{G}_{s_k} &\triangleq b -A (\affop_k H  G_0 + \affconst_k b) \\  
	\tilde{F}_{\sigma_k} &\triangleq A \affop_k H (F^*_0 - F_0) \\
	\tilde{G}_{\sigma_k} &\triangleq A \affop_k H (G^*_0 - G_0) 
  \end{align}
\end{subequations}
The following lemma makes the relationship between $\tilde\alpha^j_k$ and $\alpha^j_k$ more explicit
\begin{lemma}
  \label{prop:alp} If no constraint has been removed by Algorithm \ref{alg:primal-as} up until iteration $k$, 
  $\tilde\alpha^j_k = \tau + (1-\tau) \alpha^j_k$, $\tau \in [0, 1).$  
\end{lemma}

\begin{IEEEproof}
  Since only constraints have been added to $\mathcal{W}$ since Algorithm \ref{alg:primal-as} started, it follows from Corollary \ref{prop:add-col} and \ref{prop:add-col2} that $x_{k} = \affop_k H (x_0+\tau p_0)+\affconst_k b(\theta)$  and ${p_k = (1-\tau) \tilde H^*_k H p_0}$ for some $\tau \in [0,1)$. This inserted into \eqref{eq:step-length} gives  
  \begin{equation*}
	\begin{aligned}
	  \alpha_k^j(\theta) &= \frac{[b(\theta)]_j-[A]_j(\affop_k H x_0+\affconst_k b(\theta)) -\tau[A]_j \affop_k H p_0(\theta)}{[A]_j (1-\tau)\affop_k H p_0(\theta)} \\ &=  \frac{1}{1-\tau} \tilde\alpha^j_k(\theta) - \frac{\tau}{1-\tau}
	\end{aligned}
  \end{equation*}
  which is equivalent to $\tilde\alpha_k^j(\theta) = \tau + (1-\tau)\alpha_k^j(\theta)$.
\end{IEEEproof}
Next, we use Lemma \ref{prop:alp} to prove that $\Theta^{\text{+}j}_k$ can be equivalently expressed in terms of $\tilde{\alpha}^j_k$ instead of $\alpha^j_k$
\begin{lemma}
  \label{lem:TH-equiv}
  If no constraint has been removed in Algorithm \ref{alg:primal-as} up until iteration $k$, $\Theta^{\text{+} j}_k$ defined by \eqref{eq:th-add} equals 
  \begin{equation}
	\label{eq:theta-tilde-ex}
	\{\theta \in \Theta_k| [s^*_k]_j < 0,\: \tilde{\alpha}_k^j(\theta) < \tilde{\alpha}_k^i(\theta),\forall i \in \bar{\mathcal{W}}_k^{-}\setminus\{j\} \} 
  \end{equation}
\end{lemma}

\begin{IEEEproof}
  From Lemma \ref{prop:alp} we have that 
  \begin{equation}
    \begin{aligned}
	  \tilde{\alpha}^j_k < \tilde{\alpha}^i_k \Leftrightarrow& \tau +(1-\tau) \alpha^j_k < \tau + (1-\tau) \alpha^i_k \\
	  \Leftrightarrow&(1-\tau) \alpha^j_k < (1-\tau)\alpha^i_k \Leftrightarrow \alpha^j_k < \alpha^i_k 
    \end{aligned}
  \end{equation}
  where the last equivalence follows from $(1-\tau) > 0$ since $\tau \in [0,1)$. Hence we can replace $\alpha^j_k$ and $\alpha^i_k$ with $\tilde{\alpha}^j_k$ and $\tilde{\alpha}^i_k$, respectively, in \eqref{eq:th-add}.
\end{IEEEproof}

$\Theta^{\text{+}j}_k$ can now be explicitly stated, by inserting \eqref{eq:sur-steplength} in \eqref{eq:theta-tilde-ex} and rearranging terms to remove the fractions, as all $\theta \in \Theta_k$ satisfying 
\begin{subequations}
  \label{eq:th_add_case1}
  \begin{align}
	\label{eq:quad_const}
	\theta^T Q_k^{j,i} \theta + R_k^{j,i} \theta + S_k^{j,i} &< 0, \quad \forall i \in \bar{\mathcal{W}}_k\setminus\{ j\}   \\
	\label{eq:non-neg-constr-case1}
	[F^*_{s_k}]_j \theta + [G^*_{s_k}]_j &<0
  \end{align}
\end{subequations}
with $Q_k^{j,i}, R_k^{j,i}$ and $S_k^{j,i}$ defined as 
\begin{subequations}
  \begin{align}
	\label{eq:quad_term}
	Q_k^{j,i} \triangleq& [\tilde{F}_{\sigma_k}]_i^T [\tilde{F}_{s_k}]_j - [\tilde{F}_{\sigma_k}]_j^T [\tilde{F}_{s_k}]_i \\
	\begin{split}
	R_k^{j,i} \triangleq& [\tilde{G}_{s_k}]_j [\tilde{F}_{\sigma_k}]_i + [\tilde{G}_{\sigma_k}]_i [\tilde{F}_{s_k}]_j\\ 
			  &-([\tilde{G}_{s_k}]_i [\tilde{F}_{\sigma_k}]_j + [\tilde{G}_{\sigma_k}]_j [\tilde{F}_{s_k}]_i)
	\end{split}\\
	S_k^{j,i} \triangleq& [\tilde{G}_{s_k}]_j [\tilde{G}_{\sigma_k}]_i - [\tilde{G}_{s_k}]_i [\tilde{G}_{\sigma_k}]_j
  \end{align}
\end{subequations}
where \eqref{eq:non-neg-constr-case1} ensures that  $[A]_j x \leq [b(\theta)]_j$ is a blocking constraint and \eqref{eq:quad_const} ensures that it is the \textit{first} blocking constraint.
Thus, the parameter space will be partitioned by linear and, in contrast to Case 1, quadratic inequalities when a constraint is added to the working set under Case 2. Quadratic inequalities make the analysis less tractable compared to only linear inequalities. Hence, we will give some alternatives to circumvent these in Section \ref{ssec:spec-cases} and \ref{ssec:outer-approx}.

The results from Section \ref{sssec:case1} and \ref{sssec:case2} are summarized in Algorithm \ref{alg:add}, which describes how regions are partitioned in mode b).

%\begin{algorithm}[H]
%	\caption{Check for local optimality and add constraints to $\mathcal{A}$}\label{alg:add}
%  \begin{algorithmic}[1]
%	\State $\verb!++!N_i$
%	\State Compute $F^*$ and $G^*$ according to \eqref{eq:para-kkt}
%	\If {$\hat{p} ==$ NaN}
%	\For {all $i$ in $\mathcal{I}$ }
%	\State Calculate $\Theta^i$ according to \eqref{eq:add-part-1} 
%	\If {$\Theta^i \neq \emptyset $ } 
%	\State Append $(\Theta^i, \mathcal{A}\cup{i}, *, *, \text{f},\text{f},N_i, \hat{p})$ to $P$
%	\EndIf
%  \EndFor
%  \State Calculate $\Theta^*$ according to \eqref{eq:local-part-1} 
%	\If {$\Theta^* \neq \emptyset $ } 
%	\State Append $(\Theta^*, \mathcal{A},F^*,G^*,\text{t},\text{f},N_i, \hat{p})$ to $P$
%  \EndIf
%  \Else
%	\For {all $i$ in $\mathcal{I}$ }
%	\State Calculate $\Theta^i$ according to \eqref{eq:add-part-2} 
%	\State Calculate $F^+_i$ and $G^+_i$ according to \eqref{eq:x-2} 
%	\If {$\Theta^i \neq \emptyset $ } 
%	\State Append $(\Theta^i, \mathcal{A}\cup{i},F^+_i,G^+_i,\text{f},\text{f},N_i, \hat{p})$ to $P$
%	\EndIf
%  \EndFor
%  \State Calculate $\Theta^*$ according to \eqref{eq:local-part-2} 
%	\If {$\Theta^* \neq \emptyset $ } 
%	\State Append $(\Theta^*, \mathcal{A},F^*,G^*,\text{t},\text{f},N_i, \hat{p})$ to $P$
%  \EndIf
%  \EndIf
%  \State \Return {$P$}
%  \end{algorithmic}
%  \todo[inline]{Treat the different cases in "parallel" to save space.}
%\end{algorithm}
\begin{algorithm}[H]
  \caption{Partition $\Theta$ based on if a CSP is reached or if a constraint is added to $\mathcal{W}$. (Case 1/Case 2)}\label{alg:add}
  \begin{algorithmic}[1]
	\State \textsc{ModeB}{($(\Theta,\mathcal{W},F,G,s,k, \hat{n})$ , mpQP)}
	%\Require $(\Theta,\mathcal{A},F,G,s,k, \hat{p})$ , mpQP
	%\Ensure $P$
	\State $\verb!++!k$
	\State Compute $F^*$ and $G^*$ according to \eqref{eq:F-star}
	\For {all $i$ in $\bar{\mathcal{W}}$ }
	\State Calculate $\Theta^{\text{+}i}$ according to \eqref{eq:th_add_case2}/\eqref{eq:th_add_case1}
	\State Calculate $F_+$ and $G_+$ according to \eqref{eq:x_plus}/(-)
	\If {$\Theta^{\text{+}i} \neq \emptyset $ } \label{line:alg3-empty1} 
	\State Append $(\Theta^{\text{+}i}, \mathcal{W}\cup\{i\},F_+,G_+,0,k, \hat{n})$ to $P$
	\EndIf
  \EndFor
  \State Calculate $\Theta^{\text{CSP}}$ according to \eqref{eq:opt-slack-part} 
  \If {$\Theta^{\text{CSP}} \neq \emptyset $ } \label{line:alg3-empty2}
  \State Append $(\Theta^\text{CSP}, \mathcal{W},F^*,G^*,1,k, \hat{n})$ to $P$
  \EndIf
  \State \Return {$P$}
  \end{algorithmic}
\end{algorithm}

  \begin{remark}
  Efficient active-set solvers perform low-rank modifications to the factorization of relevant matrices when a constraint is removed or added to $\mathcal{W}$ \cite{nielsen2017low}. The same factorization techniques can be used to decrease the computational complexity of Algorithm \ref{alg:remove} and \ref{alg:add}. \end{remark}
\begin{remark}
  As was previously mentioned there are many different primal active-set algorithms in the literature and numerous of these methods are equivalent in the sense that they produce the same sequence of iterates \cite{best1984equivalence}. The main difference between algorithms is how, and which, matrices are factorized for solving the KKT-system. Hence, to determine the FLOPs for a specific algorithm one simply needs a mapping $\mathcal{F}(\mathcal{W}_k)$ that takes a working set and calculates the needed number of flops to compute the search direction. This allows for \textit{simultaneous} comparison of the FLOPs for, e.g., null-space, range-space and full-space methods when Algorithm \ref{alg:primal-as} is applied to \eqref{eq:mpQP}. Hence, the choice of, e.g., the method for solving linear equations systems can be optimized w.r.t. to the specific problem at hand.
\end{remark}
\subsection{Special cases}
\label{ssec:spec-cases}
As has been shown in \eqref{eq:th_add_case1}, the application of the proposed method to a general mpQP might result in a partitioning of the parameter space using not only affine but also quadratic inequalities. The significance of this is during the pruning of empty regions, done at line \ref{line:alg4-empty1} and \ref{line:alg4-empty2} of Algorithm \ref{alg:remove} and line \ref{line:alg3-empty1} and \ref{line:alg3-empty2} of Algorithm \ref{alg:add}, since to check consistency of a combination of linear and quadratic constraints is non-trivial. 
However, there are some relevant cases when the partitioning is solely composed of affine constraints, resulting in an easier analysis since to check whether an intersection of half-planes is empty or not can be done by solving one LP. Such special cases are described below.   
\subsubsection{No state constraints}
\label{sssec:no-state-constraint}
When there are no constraints on the states, a linear MPC problem can be formulated as an mpQP with $W=0$.  Additionally, an admissible control input can be picked as a fixed starting point, i.e., $F_0 = 0$. This will result in $[\tilde{F}_k]_j=0$ in \eqref{eq:Cj_case1} which in turn results in $Q_k^{j,i} = 0$ in \eqref{eq:quad_term}. Therefore, all partitioning of the parameter space will be done using half-planes, leading to a polytopic partition.

\subsubsection{Starting in a constrained stationary point}
\label{sssec:start-csp}
When the initial point is a constrained stationary point, partitioning according to Case 2 will never occur. Hence, under the assumption that $\Theta_0$ is a polyhedron, the final partition will be polytopic since all further partitioning of the parameter space in Case 1 is done by half-planes, see Remark \ref{rem:case1nice}.
\subsubsection{Reformulate QP using a quadratic penalty method}
All inequality constraints that depend on parameters can be transformed to equality constraints by introducing slack variables. These equality constraints can then be moved to the objective function under a quadratic penalty, cf. e.g., \cite{saraf2018fast}\cite[Sec.17-1]{nocedal}.
The resulting QP will be on the form which was discussed in Case 1, described above. 

\subsection{Outer approximations of quadratic inequalities}
\label{ssec:outer-approx}
The comparison of step lengths $\alpha^i_k < \alpha^j_k$ to find the first blocking constraint when a constraint is yet to be removed from $\mathcal{W}$, i.e. under Case 2, results in the quadratic inequalities \eqref{eq:quad_const} on the form  
\begin{equation}
  \theta^T Q \theta + R \theta + S < 0  
\end{equation}
  As previously mentioned, the consistency check that is done in Algorithm \ref{alg:remove} and \ref{alg:add} will be more challenging when both affine and quadratic inequalities define a region, in particular since $Q$ can be indefinite. An alternative to these quadratic constraints is to make an affine outer-approximation with the half-plane  
\begin{equation}
  \label{eq:relax-quad}
  R \theta <-S - \min_{\theta\in\Theta_k} \theta^T Q \theta 
\end{equation}
where $\Theta_k$ is the current region. Hence, by solving an indefinite QP in relatively low dimension, an affine relaxation can be obtained. Ultimately, relaxing the quadratic constraints might lead to some regions overlapping, giving a conservative result since all regions produced by the certification method might not correspond to how the Algorithm \ref{alg:primal-as} performs in practice. 

An interpretation of  relaxing $\alpha^i_k < \alpha^j_k$ with \eqref{eq:relax-quad} is that the $i$:th constraint might not be the \textit{first} blocking constraint for that particular parameter region in iteration $k$. This would result in an primal infeasible iterate, which can be used in the certification algorithm to prune some of the redundant regions which the outer-approximation might yield. Checking the infeasibility of the iterate during Case 2 will, again, lead to quadratic regions and is therefore of no use. However, as soon as a constrained stationary point is reached, the iterates become affine in $\theta$, see Theorem \ref{prop:affine},  and the affine constraints $A x^*_k (\theta) \leq b + W\theta$ can be added to the current region to prune infeasible iterates. In the end, the only redundant regions that remain will correspond to iterates that regained primal feasibility before the first CSP was reached. 
\section{Extension to positive semi-definite case}
\label{sec:singular}
We will now extend Algorithm \ref{alg:primal-as} and its corresponding certification method Algorithm
\ref{alg:main} to the positive semi-definite case, i.e., $H \succeq 0$. Not only does this allow the primal active-set algorithm to be certified for a more general problem class, it also allows us to certify a family of dual active-set methods and active-set methods used in linear programming, creating a unifying framework for certification of active-set methods which change a single index at a time in the working set.   

\subsection{Extending the active-set algorithm}
\label{ssec:ext-alg}
If the Hessian of the QP being solved is positive semi-definite, the reduced Hessian $Z_k^T H Z_k$ can become singular. In that case, $p_k$ cannot be computed by simply solving the KKT-system \eqref{eq:kkt-step} and must be determined in another way \cite{gill2015methods}.  $Z_k^T H Z_k$ being singular means that the objective function on the subspace defined by the current working set lacks a quadratic part, i.e., is affine and therefore is in general unbounded along a direction in the subspace. Such a direction can be found by solving the system 

\begin{equation}
  \label{eq:kkt-singular}
  \begin{pmatrix}
	H & A^T_k \\ A_k^T & 0 
  \end{pmatrix}
  \begin{pmatrix}
	p_k \\ \lambda_k  
  \end{pmatrix}
  = 
\begin{pmatrix}
 0 \\ 0  
\end{pmatrix}
\end{equation}
and an example of a solution to this KKT-system is given by the following lemma 
\begin{lemma}
  \label{lem:sing-direction}
  If $Z_k^T H Z_k$ becomes singular after removing the $i$:th row of $A_{k-1}$, a solution $p_k$ to \eqref{eq:kkt-singular} is given by $T_{k-1} e_i$, where $T_k$ is defined in \eqref{eq:Hstar-def} and $e_i$ is the $i$:th unit vector. 
\end{lemma}
\begin{IEEEproof}
  Let the $i$:th row of $A_{k-1}$, which is removed, be denoted $\tilde{a}^T$. Furthermore, let $\Pi$  be a permutation matrix for which  $\Pi A_{k-1} = [A_k^T, \:  \tilde{a}]^T$, i.e., a permutation matrix which moves the $i$:th row to the last row. Then we have that  
  \begin{equation*}
	\begin{pmatrix}
	  I & 0 \\ 0 & \Pi
	\end{pmatrix}
	\begin{pmatrix}
	  H & A_{k-1}^T \\ A_{k-1} & 0 
	\end{pmatrix}
	\begin{pmatrix}
	  I & 0 \\ 0 & \Pi^T
	\end{pmatrix}
	=
    \begin{pmatrix}
	  H & A_{k}^T & \tilde{a} \\
	  A_{k} & 0 & 0 \\
	  \tilde{a}^T & 0 & 0 \\ 
    \end{pmatrix}
  \end{equation*}
  which is nonsingular since the KKT-system at iteration $k-1$ in nonsingular. Taking the inverse of this matrix gives  
  \begin{equation}
	\label{eq:perm-inv}
    \begin{pmatrix}
	  H & A_{k}^T & \tilde{a} \\
	  A_{k} & 0 & 0 \\
	  \tilde{a}^T & 0 & 0 \\ 
	\end{pmatrix}^{-1} =  
   \begin{pmatrix}
	 H^*_{k-1} & T_{k-1} \Pi^T \\
	 \Pi T^T_{k-1}  & \Pi U_{k-1} \Pi^T 
   \end{pmatrix} 
  \end{equation}
  where we have used \eqref{eq:kkt-inv} and $\Pi^T = \Pi^{-1}$. Now, consider the system 
  \begin{equation}
	\label{eq:kkt-ext}
    \begin{pmatrix}
	  H & A_{k}^T & \tilde{a} \\
	  A_{k} & 0 & 0 \\
	  \tilde{a}^T & 0 & 0 \\ 
    \end{pmatrix}
	\begin{pmatrix}
	  p_k \\ \lambda_k \\  0  
	\end{pmatrix}
	= 
	\begin{pmatrix}
	0 \\ 0 \\ 1 
	\end{pmatrix}
  \end{equation}
  and note that the two first rows is equivalent to \eqref{eq:kkt-singular}. As a result, a solution to \eqref{eq:kkt-ext} is a solution to \eqref{eq:kkt-singular}. Multiplying both sides of \eqref{eq:kkt-ext} from left with \eqref{eq:perm-inv} gives 
  \begin{equation}
	p_k = T_{k-1} \Pi^T \begin{pmatrix}
	 0 \\ 1  
	\end{pmatrix}
	= T_{k-1} e_i 
  \end{equation}
\end{IEEEproof}

When deciding the step length for the singular case, two different scenarios can occur. If there is a blocking constraint along the ray $x_k + \alpha p_k, \alpha> 0$, the blocking constraint can be added to the working set, and the iterations can proceed as usual. Otherwise, if there are no blocking constraints along the ray, the objective function can be decreased by an arbitrary amount by moving along the ray, since it is a descent direction, resulting in an unbounded problem. Concretely, there will be no blocking constraint if $A p_k \triangleq \sigma_k \geq 0$ since then the updated slack $s_{k+1} = s_k + \alpha \sigma_k$ cannot become negative for any positive step length $\alpha$, i.e., any positive $\alpha$ gives a primal feasible iterate. 
The modifications described above are summarized in Algorithm \ref{alg:AS-sd-iter}.

If $Z_0^T H Z_0$ is nonsingular, $Z_k^T H Z_k$ will have at most one singular eigenvalue by the following reasoning. If it is singular in iteration $k$, the scheme outlined above either adds a constraint to $\mathcal{W}$, if a blocking constraint exists, which cannot introduce more singular eigenvalues to the reduced Hessian. Otherwise, if there are no blocking constraints, the problem is marked as unbounded and the algorithm terminates. 
\begin{algorithm}  
  \caption{Iteration in Algorithm \ref{alg:primal-as} when $Z_k^T H Z_k$ is singular}
  \label{alg:AS-sd-iter}
  \begin{algorithmic}[1]
	\State Compute $p_k$ from \eqref{eq:kkt-singular} 
	\State $[\sigma_k]_{\mathcal{C}_k} \leftarrow [A]_{\mathcal{C}_k} p_k$
	\If {$\sigma_k \geq 0$} 
	\State \textbf{break} unbounded
	\Else\quad $m \leftarrow \underset{i \in \mathcal{C}_k:[\sigma_k]_i > 0}{\text{argmin }} \frac{[s_k]_i}{[\sigma_k]_i};\:\:\: \mathcal{W}_{k+1}  \leftarrow \mathcal{W}_k \cup \{m\}$  
	\State \quad $x_{k+1} \leftarrow x_k + \alpha^m_k p_k$;\:\:\:\:\: $s_{k+1} \leftarrow s_k - \alpha^m_k \sigma_k$
	\EndIf
	\State $k \leftarrow k+1$ 
  \end{algorithmic}
\end{algorithm}

\begin{remark}
The method employed when the reduced Hessian is singular can be seen as a switching rule for the working set. The reduced Hessian becomes singular after a constraint has been removed, and this will always lead to another constraint being added, assuming that the problem is bounded, which can be seen as a "switch" of indices in the working set. 
\end{remark}
\subsection{Extending the certification algorithm}  
Since a normal iteration of Algorithm \ref{alg:primal-as} can be performed when $Z_k^T H Z_k$ is nonsingular, amendments to the certification algorithm only need to be considered when $Z_k^T H Z_k$ is singular. Moreover, since $Z_k^T H Z_k$ only becomes singular after a constraint has been removed, modifications only have to be made for mode b).  

In the singular case, $p_k$ is independent of the parameter since it is computed by solving \eqref{eq:kkt-singular} which does not contain $\theta$. As was discussed in the previous section, if $ A p_k \triangleq \sigma_k \geq 0$ there are no blocking constraints, resulting in an unbounded problem, hence, we mark the region $\Theta_k$ as unbounded if $\sigma_k$ has no negative components. Otherwise we will have blocking constraints, corresponding to the negative components, and for these we partition the parameter space depending on the \textit{first} blocking constraint. Explicitly, the region for which the $j$:th constraint is the first blocking constraint, and hence will be added to $\mathcal{W}$, is

\begin{equation}
  \begin{aligned}
	\Theta^j_k =& \{\theta \in \Theta_k |\alpha_k^j(\theta) < \alpha_k^i(\theta),\:\: \forall i : [\sigma_k]_i  < 0\} \\
  \end{aligned}
\end{equation}
which, analogously to what was described in Section \ref{ssec:add}, can be written as all $\theta \in \Theta_k$ such that  
\begin{equation}
 \begin{aligned}
	&K^{j,i}_k \theta <L^{j,i}_k, \quad \forall i: [\sigma_k]_i < 0  \\
 \end{aligned} 
\end{equation}
with the same definitions of $K$ and $L$ as in \eqref{eq:KL-def} except that $G_{\sigma_k} \triangleq A p_k$.
\subsection{Dual active-set methods for Quadratic Programming}
With the extension to semi-definite problems, we now turn our attention to dual active-set QP methods. 
As is noted in \cite[p.244]{fletcher} and \cite{10.1007BF02591962}, the popular dual active-set method presented in \cite{10.1007BF02591962}, which we will call the Goldfarb-Idnani (GI) method, is equivalent to Algorithm \ref{alg:primal-as}, with the extensions mentioned in Section \ref{ssec:ext-alg},  being applied to the dual of \eqref{eq:mpQP} when $H \succ 0$. The dual problem to \eqref{eq:mpQP} can be stated as the following mpQP 
\begin{equation}
  \label{eq:qp_dual}
  \begin{aligned}
	 &\underset{\lambda}{\text{minimize}}&& \frac{1}{2}\lambda^T A^T H^{-1} A \lambda + (f^T(\theta)H^{-1} A^T +b^T(\theta)) \lambda\\
     &\text{subject to} && \lambda \geq 0.
  \end{aligned}
\end{equation}
where the optimal primal solution $x^*$ is related to, and can be recovered from, the optimal dual solution $\lambda^*$ by ${x^* = -H^{-1}(f(\theta) +A^T \lambda^*)}$.

A complexity certification method for the GI method is provided in \cite{10.1109TAC.2017.2696742}, where the number of iterations is shown to be constant over a polyhedral partition of the parameter space. This is in contrast with the results in Section \ref{sec:cert-primal-as} where both affine \textit{and} quadratic inequalities partitions the parameter space for Algorithm \ref{alg:main}.
There are two factors that, separately, lead to a partition solely of polyhedral type for the dual active-set method. First, in \cite{10.1109TAC.2017.2696742} the dual active-set method is always initialized in the unconstrained optimum, which implies that all dual variables are $0$ in the first iteration and all constraints are active, i.e., the first iterate is a constrained stationary point. This falls into the special case discussed in Section \ref{sssec:start-csp}, which results in a polyhedral partition.  

A second reason for a final polyhedral partition is that \eqref{eq:qp_dual} has more structure than the generic mpQP in \eqref{eq:mpQP}, namely that there is no parameter dependence in the constraints. This additional structure will, with the same reasoning as in the special case described in Section \ref{sssec:no-state-constraint}, lead to a polyhedral partition, even if the method is not started in a constrained stationary point (as long as this starting point is parameter independent). 

The certification of a dual active-set method that is not started in the unconstrained optimum is not considered in \cite{10.1109TAC.2017.2696742}. However, viewing the method as Algorithm \ref{alg:main}, with the amendments to handle the singular case, applied to the dual makes it possible to certify a dual active-set method that starts with an arbitrary, dual feasible, starting iterate. Being able to do the certification from an arbitrary starting iterate is necessary when analyzing the behaviour of the method when it is warm-started. 

\subsection{Active-set methods for Linear Programming}
Using another formulation, more concretely using the $1$- and $\infty$-norm instead of the $2$-norm in the cost function, linear MPC problems can be cast as mpLPs, see, e.g.,  \cite[Sec.2-3]{borrelli2003constrained}.  
mpLPs can be seen as a special class of mpQPs with $H=0$. A well-renowned method for solving LPs is the simplex method \cite[Sec. 5]{dantzig1963linear} which is also an active-set method. In fact, Algorithm \ref{alg:primal-as}, with the amendments from Section \ref{ssec:ext-alg}, applied to an LP is equivalent to the simplex method with Dantzig's pivot rule \cite{gill2015methods}, where equivalent means that the same iterate sequences are produced by both methods. The iterates of the simplex method are vertices of the feasible set and we will now briefly describe how this translates to the behaviour of Algorithm \ref{alg:primal-as} with its singular extension. Since a vertex is a CSP, we will check for optimality or remove a constraint from our working set (mode a)). Removing a constraint leads to a singular reduced Hessian which, in turn,  leads to a computation of the step direction according to \eqref{eq:kkt-singular}. As was discussed before, this search will either lead to no constraint being encountered along $p_k$, in which case the problem is unbounded, or a constraint will be encountered and added to the working set, resulting in a new vertex. 

As an alternative to the simplex method for solving LPs, one can use another active-set algorithm which does not restrict all iterates to vertices. Such a method is considered in \cite{10.1109TAC.2011.2108450} and uses the gradient of the objective function as search direction. Using the gradient as a search direction results in the KKT-system  
\begin{equation}
  \label{eq:LP-step}
  \begin{pmatrix}
	I & A^T_k \\ A_k & 0 
  \end{pmatrix}
  \begin{pmatrix}
   p_k \\ \lambda_k 
\end{pmatrix} = 
\begin{pmatrix}
  -f(\theta)\\
  0
\end{pmatrix}.
\end{equation}
Computing $p_k$ by \eqref{eq:LP-step} instead of \eqref{eq:kkt-singular} in Algorithm \ref{alg:AS-sd-iter} leads to this LP algorithm. In \cite{10.1109TAC.2011.2108450}, this active-set method was certified for mpLPs with $f_{\theta} = 0$, i.e., $f(\theta) = f$. 
\section{NUMERICAL EXAMPLES}
\label{sec:examples}
Some benchmark problems from the MATLAB Model Predictive Control Toolbox were considered to test the proposed certification method. These MPC problems were the control of a double integrator, a DC-motor, an inverted pendulum, a linearized nonlinear multiple-input-multiple-output system and an ATFI-F16 aircraft. The tracking problem was considered, resulting in a parameter vector $\theta$ containing the state vector, the previous control input and the reference signal. The same problems were also considered in the context of real-time certification for other QP methods in \cite{10.1109TAC.2017.2696742} and \cite{cimini2019complexity}, where they were considered a good representation of the kind of problems encountered in real-time MPC. For further details about the problems see \cite{10.1109TAC.2017.2696742} and \cite{cimini2019complexity}. Additionally, the method was tested on a randomly generated mpQP to accentuate the possibility of quadratic partitioning of the parameter space. This problem is given by   

\begin{equation*}
  \begin{aligned}
	&H = \begin{pmatrix}
	0.97 & 0.19 & 0.15 \\
	0.19 & 0.98 & 0.05 \\
	0.15 & 0.05 & 0.99
  \end{pmatrix},
	&&A = \begin{pmatrix}
	0.38 & 2.20 & 0.43\\
	0.49 & 0.57 & 0.22 \\
	0.77 & 0.46 & 0.41
  \end{pmatrix}\\
	&f = \begin{pmatrix}
	  0&0&0
\end{pmatrix}^T,
	&&b = \begin{pmatrix}
	  4.1&3.7&4.3
\end{pmatrix}^T\\
	&W = \begin{pmatrix}
  0.19&-0.89\\
  0.62&-1.54\\
  -0.59&-1.01
\end{pmatrix},
	&&f_{\theta} = \begin{pmatrix}
  11.3&-44.3\\
  -3.66&-11.9\\
  -32.6&7.81\\
\end{pmatrix}
  \end{aligned}
\end{equation*}
and will be called "Contrived mpQP". 

The certification method presented in Section \ref{sec:cert-primal-as} was applied to the resulting primal mpQP problems on the form \eqref{eq:mpQP} for all of the MPC examples, with the starting iterate being the origin, i.e. $x_0 = (0, ..., 0)^T$ and the starting working set being the empty set, i.e. $\mathcal{W}_0 = \emptyset$. 
Since the DC motor and ATFI-F16 aircraft examples contain state constraint, these constraints were soften, cf. e.g., \cite{zheng1995stability}, to ensure the existence of primal feasible solutions. Furthermore, the initial slack was set large enough to ensure primal feasibility of the origin for all parameters of interest. 

In addition to the primal problems, the certification method was applied to the dual problems on the form \eqref{eq:qp_dual}, which are positive semi-definite. Hence, the amendments to the certification method described in Section \ref{sec:singular} were used. For all of the examples, the starting iterate was chosen as ${\lambda_0 = (0, ..., 0)^T}$ and all constraints of the dual problem were active in the initial working set, i.e., $\mathcal{W}_0 = \mathcal{K} = \{1,\dots,m\}$. 

Gurobi 9.0 \cite{gurobi} was used to decide if regions described by both linear and quadratic inequalities were empty or not.

\subsection{Complexity certification}
\begin{figure}
  \centering
  \input{fig/iter-bar.tex}
  \begin{subfigure}[b]{0.23\textwidth}
	\begin{tikzpicture}
	  \node[anchor=south west,inner sep=0] (image) at (0,0) {\includegraphics[width=0.95\textwidth]{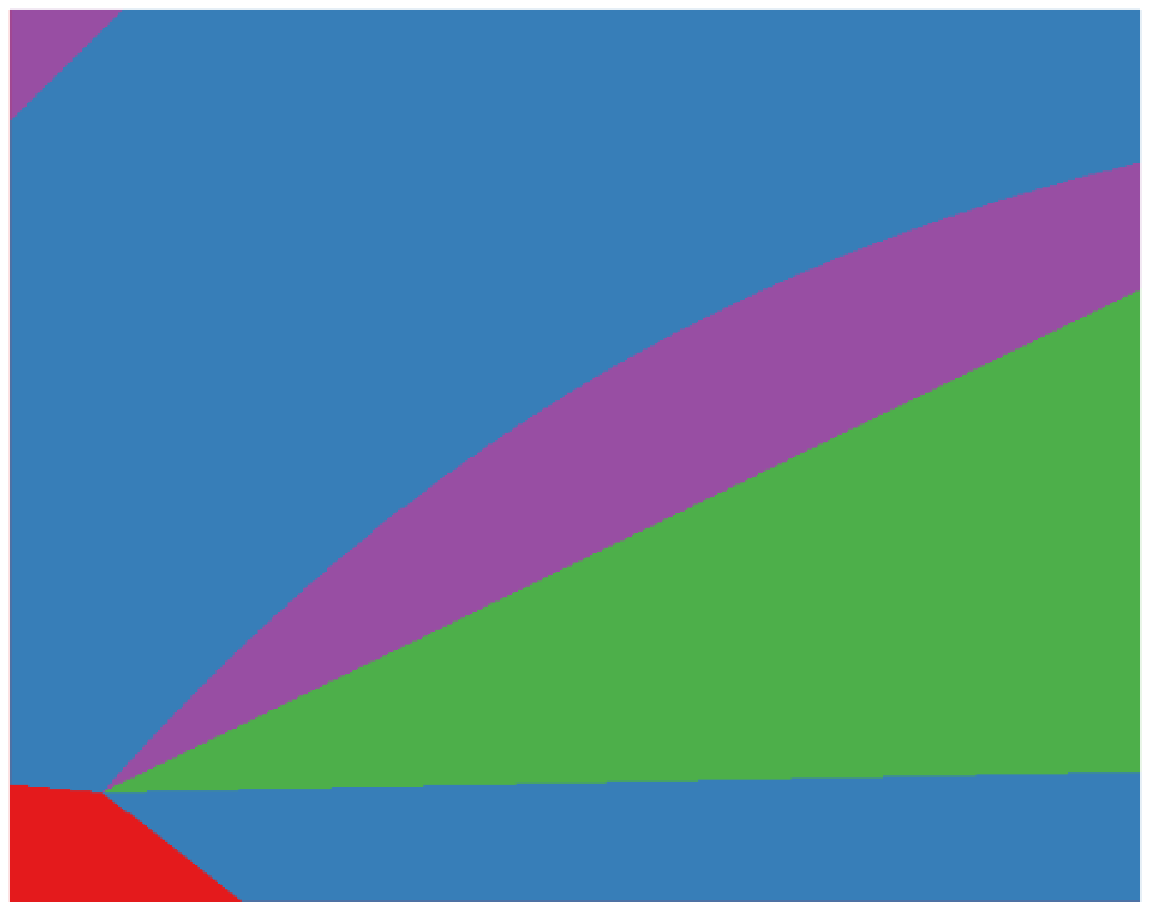}};
	  \begin{scope}[x={(image.south east)},y={(image.north west)}]
	  \node (xlabel) at (0.5,0.05) {$\theta_1$};
	  \node (xmintick) at (0.155,0.05) {$0$};
	  \node (xmaxtick) at (0.875,0.06) {$1$};
	  \node (ylabel) at (0.085,0.5) {$\theta_2$};
	  \node[align=right] (ymintick) at (0.09,0.15) {$0$};
	  \node[align=right] (ymaxtick) at (0.09,0.85) {$1$};
	  \node (dummy) at(-0.025,0.5){};
	\end{scope}
  \end{tikzpicture}
	\caption{Contrived mpQP, $p=2$}
	\label{fig:genmpqp}
  \end{subfigure}
  \begin{subfigure}[b]{0.23\textwidth}
	\begin{tikzpicture}
	  \node[anchor=south west,inner sep=0] (image) at (0,0) {\includegraphics[width=0.95\textwidth]{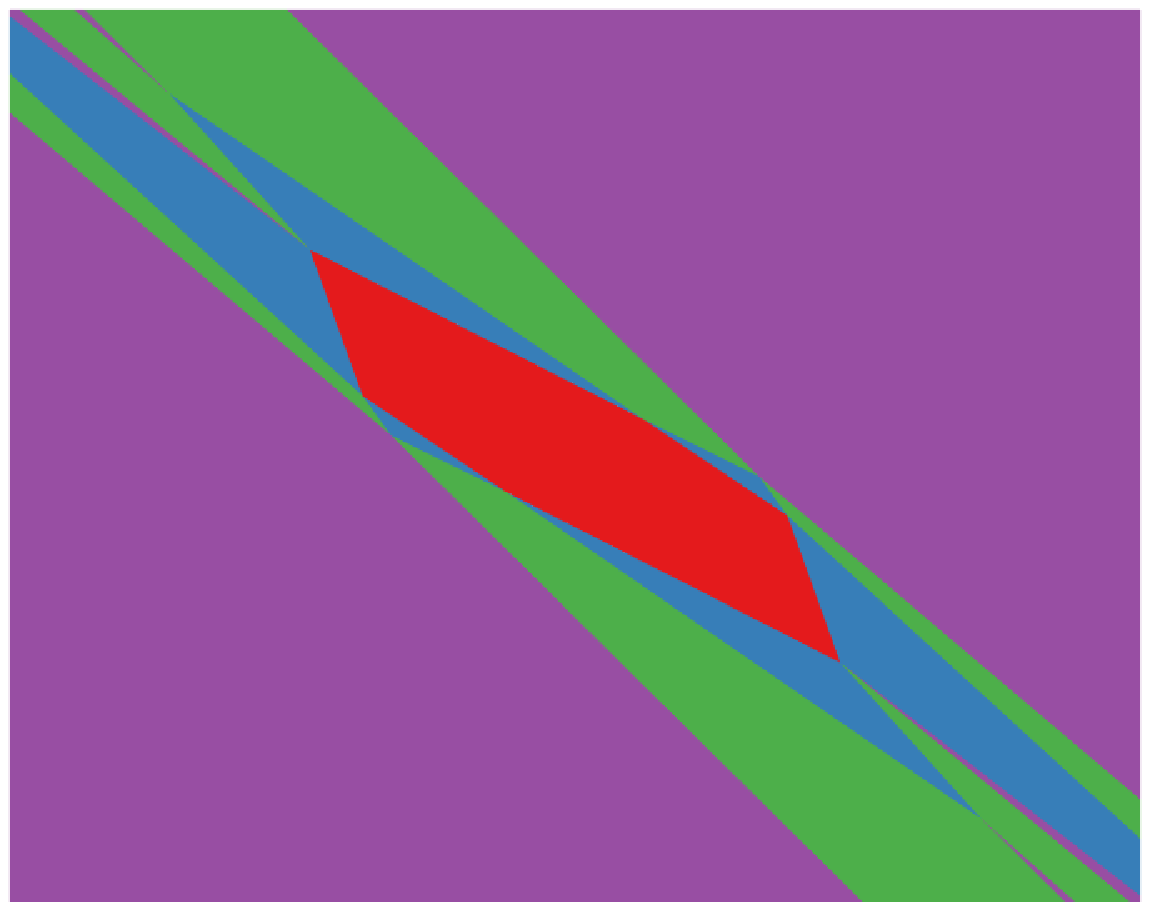}};
	  \begin{scope}[x={(image.south east)},y={(image.north west)}]
	  \node (xlabel) at (0.5,0.05) {$\theta_1$};
	  \node (xmintick) at (0.155,0.05) {$\text{-}1.5$};
	  \node (xmaxtick) at (0.875,0.06) {$1.5$};
	  \node (ylabel) at (0.085,0.5) {$\theta_2$};
	  \node[align=right] (ymintick) at (0.075,0.15) {$\text{-}1$};
	  \node[align=right] (ymaxtick) at (0.09,0.85) {$1$};
	  \node (dummy) at(-0.025,0.5){};
	\end{scope}
  \end{tikzpicture}
	\caption{Double integrator, $p=4$}
	\label{fig:dbint}
  \end{subfigure}
  \begin{subfigure}[b]{0.23\textwidth}
	\begin{tikzpicture}
	  \node[anchor=south west,inner sep=0] (image) at (0,0) {\includegraphics[width=0.95\textwidth]{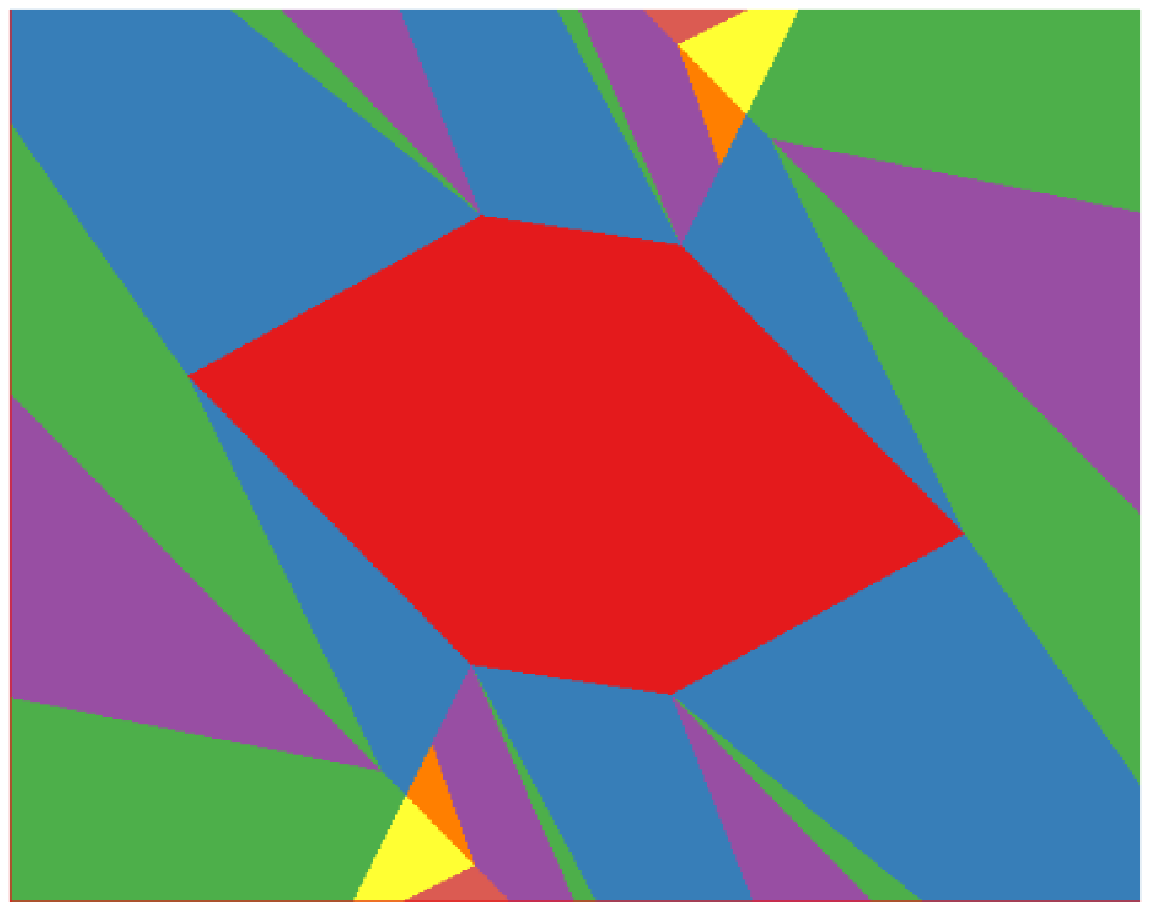}};
	  \begin{scope}[x={(image.south east)},y={(image.north west)}]
	  \node (xlabel) at (0.5,0.05) {$\theta_1$};
	  \node (xmintick) at (0.125,0.05) {$\text{-}0.15$};
	  \node (xmaxtick) at (0.875,0.06) {$0.15$};
	  \node (ylabel) at (0.085,0.5) {$\theta_2$};
	  \node[align=right] (ymintick) at (0.08,0.15) {$\text{-}1$};
	  \node[align=right] (ymaxtick) at (0.095,0.85) {$1$};
	  \node (dummy) at(-0.025,0.5){};
	\end{scope}
  \end{tikzpicture}
  \caption{DC motor, $p=6$}
	\label{fig:dcmotor}
  \end{subfigure}
  \begin{subfigure}[b]{0.23\textwidth}
	\begin{tikzpicture}
	  \node[anchor=south west,inner sep=0] (image) at (0,0) {\includegraphics[width=0.95\textwidth]{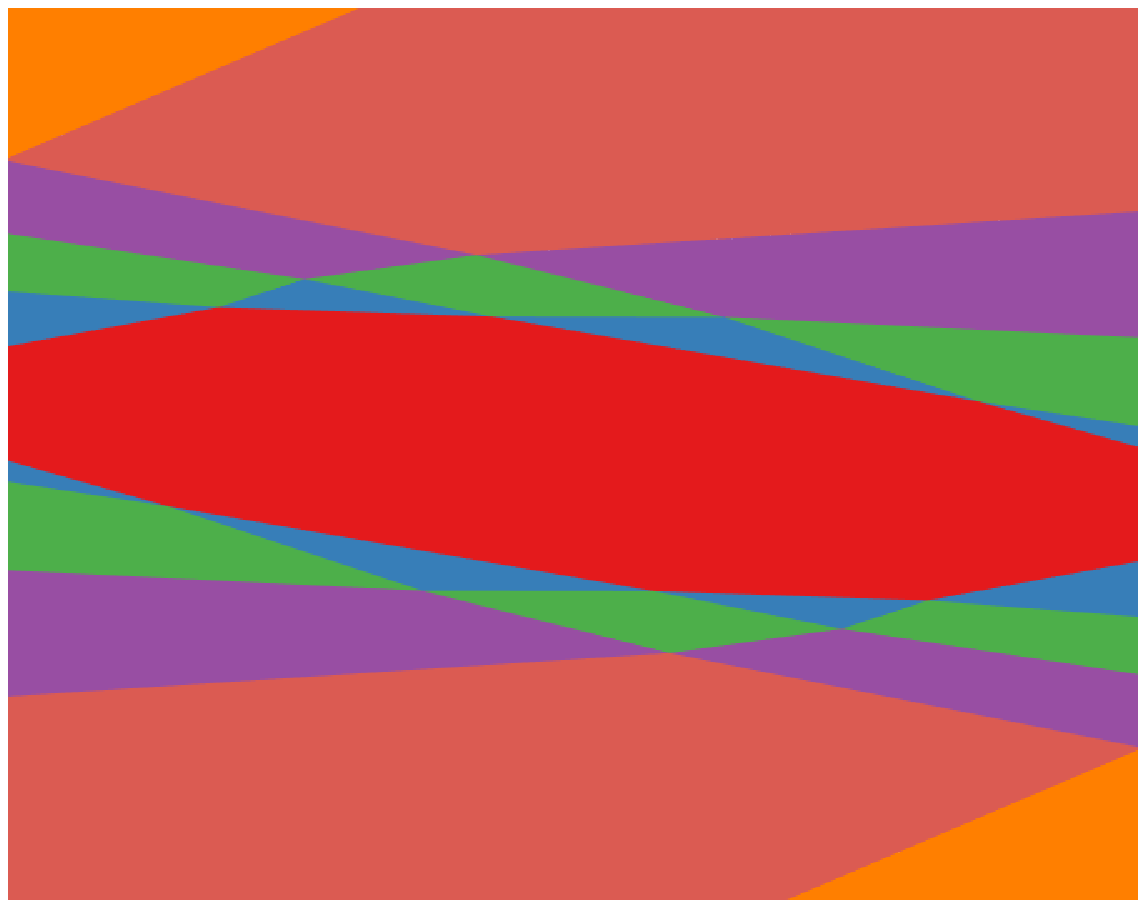}};
	  \begin{scope}[x={(image.south east)},y={(image.north west)}]
	  \node (xlabel) at (0.5,0.05) {$\theta_1$};
	  \node (xmintick) at (0.125,0.05) {$\text{-}20$};
	  \node (xmaxtick) at (0.875,0.06) {$20$};
	  \node (ylabel) at (0.085,0.5) {$\theta_2$};
	  \node[align=right] (ymintick) at (0.045,0.15) {$\text{-}20$};
	  \node[align=right] (ymaxtick) at (0.07,0.85) {$20$};
	  \node (dummy) at(-0.025,0.5){};
	\end{scope}
  \end{tikzpicture}
  \caption{Inverted pendulum, $p=8$}
	\label{fig:invpend}
  \end{subfigure}
  \begin{subfigure}[b]{0.23\textwidth}
	\begin{tikzpicture}
	  \node[anchor=south west,inner sep=0] (image) at (0,0) {\includegraphics[width=0.95\textwidth]{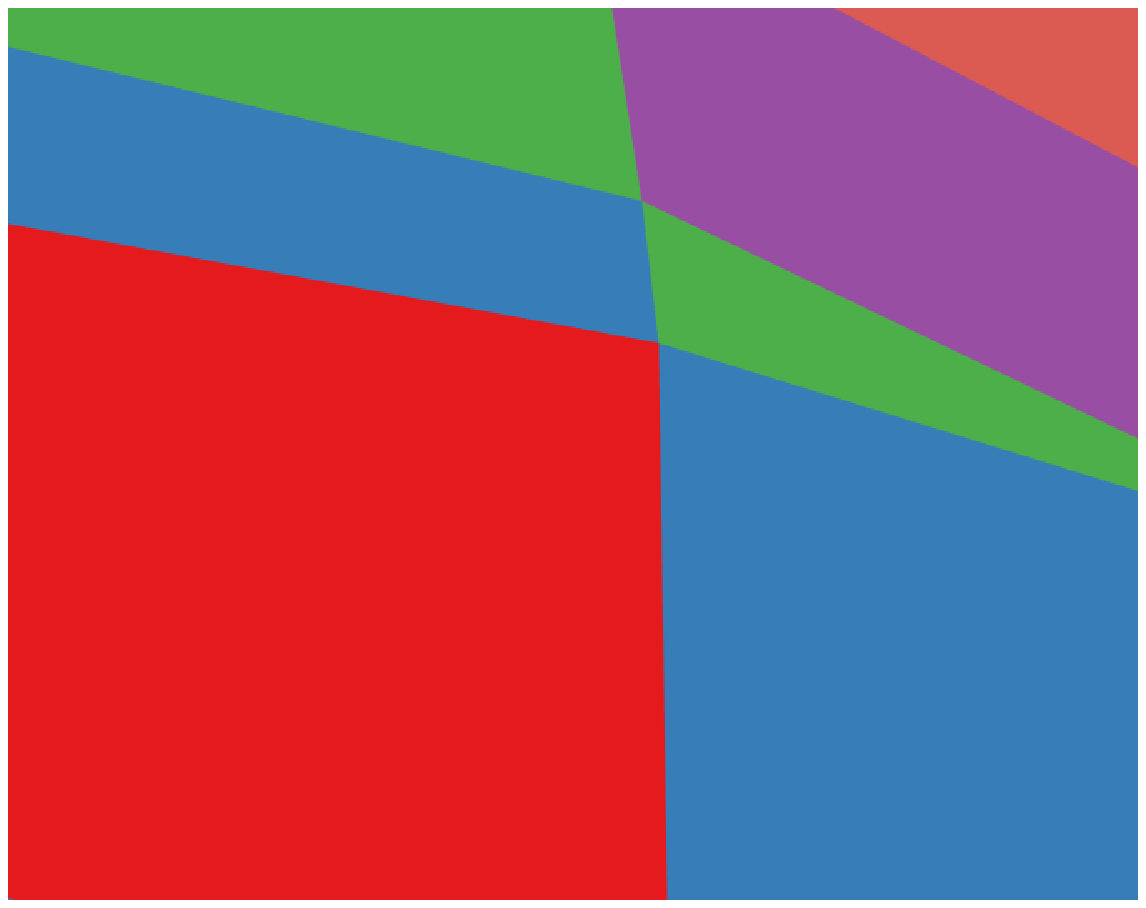}};
	  \begin{scope}[x={(image.south east)},y={(image.north west)}]
	  \node (xlabel) at (0.5,0.05) {$\theta_1$};
	  \node (xmintick) at (0.125,0.05) {$\text{-}0.5$};
	  \node (xmaxtick) at (0.875,0.06) {$2$};
	  \node (ylabel) at (0.085,0.5) {$\theta_2$};
	  \node[align=right] (ymintick) at (0.045,0.15) {$\text{-}0.5$};
	  \node[align=right] (ymaxtick) at (0.095,0.85) {$1$};
	  \node (dummy) at(-0.025,0.5){};
	\end{scope}
  \end{tikzpicture}
  \caption{Nonlinear demo, $p=10$}
	\label{fig:nonlin}
  \end{subfigure}
  \begin{subfigure}[b]{0.23\textwidth}
	\begin{tikzpicture}
	  \node[anchor=south west,inner sep=0] (image) at (0,0) {\includegraphics[width=0.95\textwidth]{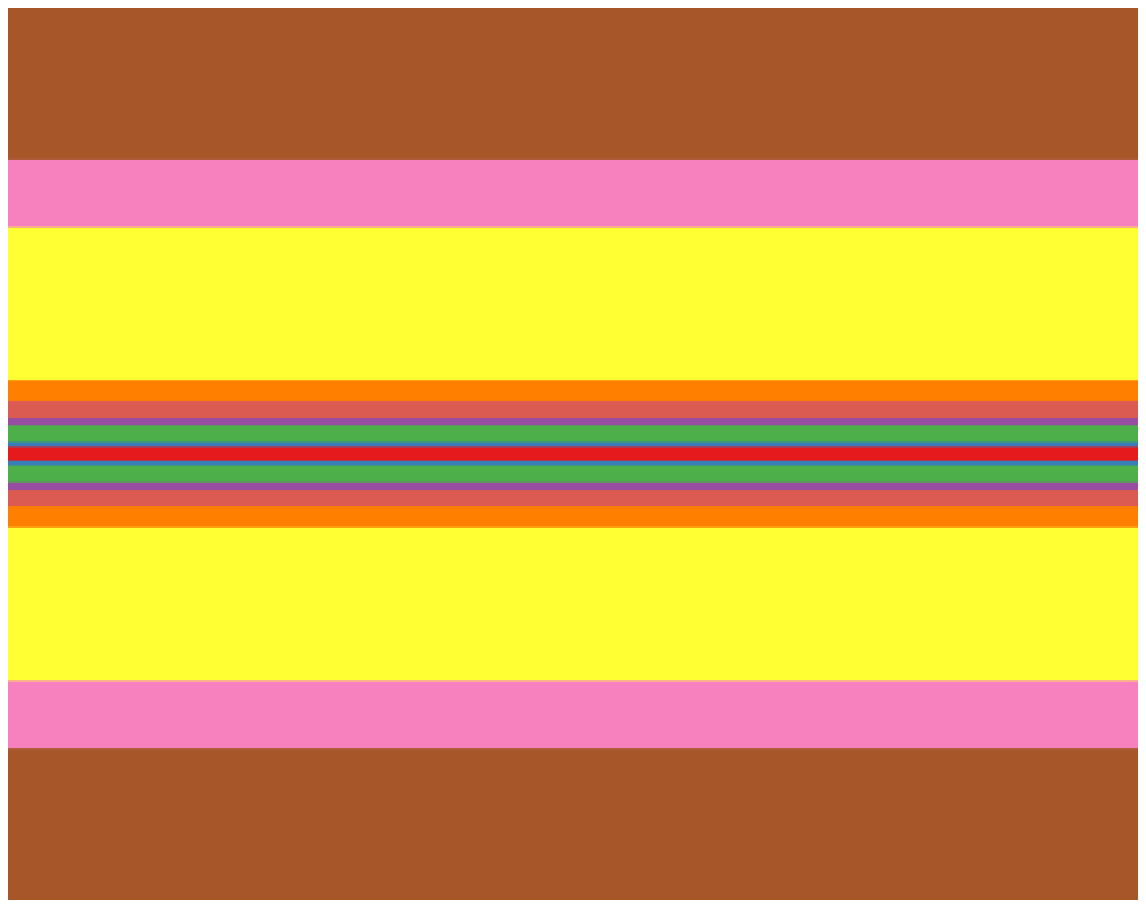}};
	  \begin{scope}[x={(image.south east)},y={(image.north west)}]
	  \node (xlabel) at (0.5,0.05) {$\theta_1$};
	  \node (xmintick) at (0.125,0.05) {$\text{-}20$};
	  \node (xmaxtick) at (0.875,0.06) {$20$};
	  \node (ylabel) at (0.085,0.5) {$\theta_2$};
	  \node[align=right] (ymintick) at (0.07,0.15) {$\text{-}20$};
	  \node[align=right] (ymaxtick) at (0.07,0.85) {$20$};
	  \node (dummy) at(-0.025,0.5){};
	\end{scope}
  \end{tikzpicture}
  \caption{ATFI-16, $p=10$}
	\label{fig:ATFI-16}
  \end{subfigure}
  \caption{2D-slice of the resulting parameter regions with ${\theta_i = 0,\: i>2}$, produced by Algorithm \ref{alg:main} for the primal problems. The same color means same number of QP iterations.}\label{fig:iterfig}

\end{figure}

To give a taste  of the final result from Algorithm \ref{alg:main}, Figure \ref{fig:iterfig} depicts a low-dimensional slice of the resulting regions which lead to the same number of QP iterations when the primal problems are solved with Algorithm \ref{alg:primal-as}, determined by Algorithm \ref{alg:main}. However, this is only a subset of the information contained in the final partition since every region also contains the exact sequence of working-set changes performed to reach the solution. As an example, the parameters in the final region of the contrived mpQP example which contains $\theta = [0.5,0.5]^T$, (the purple region in the middle of Figure \ref{fig:genmpqp}), have undergone the following working-set changes: 
${\emptyset \rightarrow \{1\} \rightarrow \{1,3\} \rightarrow \{3\}}$ before reaching optimality.

\begin{table*}[t]
  \centering
  \begin{tabular}{c||c|c|c||c|c||c|c||c|c||c|c|c}
	&$p$&$n$&$m$&$N^{\text{max}}_{\text{prim}}$&$N^{\text{max}}_{\text{dual}}$&$t^{\text{cert}}_{\text{prim}}[s]$&$t^{\text{cert}}_{\text{dual}}[s]$& $N^{\text{reg}}_{\text{prim}}$& $N^{\text{reg}}_{\text{dual}}$&$N^{\text{max}}_{\text{MC,prim}}$&$N^{\text{max}}_{\text{MC,dual}}$\\
	\hline
	Contrived mpQP&2&3&3&4&4&0.08&0.01&6&5&4&4\\
	\hline
	Double integrator&4&3&6&6&6&0.13&0.08&39&43&6&6\\
	\hline
	Inverted pendulum&8&5&10&19&14&15&7.6&2499&1839&19&14\\
	\hline
	DC motor$^*$&6&3&10&14&14&46&11&2309&1865&\textcolor{red}{10}&\textcolor{red}{10}\\
	\hline
	Nonlinear demo&10&6&12&14&11&56&41&10166&8669&\textcolor{red}{12}&11\\
	\hline
	ATFI-F16$^*$&10&5&12&21&24&541&558&41971&93064&\textcolor{red}{14}&\textcolor{red}{15}\\
	\hline
  \end{tabular}
  \newline  
  \newline  
  \footnotesize{$^*$ For the primal problem, quadratic inequalities were outer-approximated by affine inequalities as described in Section \ref{ssec:outer-approx}.}\\
  \caption{Dimensions of the resulting mpQPs for the examples, the worst-case number of QP-iterations $N^{\text{max}}$ determined by Algorithm \ref{alg:main} and the worst-case number of QP-iterations  $N^{\text{max}}_{\text{MC}}$ determined by extensive simulation. $N^{\text{reg}}$ is the number of regions in the final partition and $t^{\text{cert}}$ is the time taken by a MATLAB implementation of Algorithm \ref{alg:main} executed on an Intel$^{\tiny{\textregistered}}$ 2.7 GHz i7-7500U CPU. The subscripts "prim" or "dual" denote results when the primal or the dual QP were solved, respectively.}
  \label{tab:res}
\end{table*}%
The dimensions of the resulting mpQPs for the examples are shown in Table \ref{tab:res} together with the maximum number of QP iterations $N^{\text{max}}_{\text{primal}}$ and $N^{\text{max}}_{\text{dual}}$ needed for the active-set algorithm to provide a solution when solving the primal and dual problem, respectively, determined by Algorithm \ref{alg:main}. The table also includes the time taken for the certification $t^{\text{cert}}$ and the number of regions $N^\text{reg}$ in the final partition. Furthermore, the maximum number of QP iterations observed when running Monte Carlo (MC) simulations, denoted  $N^{\text{max}}_{\text{MC}}$, were obtained by random sampling of $\Theta_0$ and applying Algorithm \ref{alg:primal-as} to the resulting QPs. For the MC simulations, as many samples as possible were drawn during $t^{\text{cert}}$ to compare with the certification method.  

By comparing $N^{\text{max}}_{\text{prim}}$ with $N^{\text{max}}_{\text{dual}}$ in Table \ref{tab:res} it can be seen that the dual method needs fewer iterations in the worst-case for most of the examples, which is in accordance with what is noted in \cite{10.1007BF02591962}. However, for the ATFI-F16 example the primal method needs fewer iterations in the worst-case. Hence, whether the primal or dual active-set approach is to be preferred, from a real-time perspective, is, not surprisingly, problem dependent and the proposed certification method can be used to decide which one gives the fewest iterations in the worst-case for a given problem.  

It can also be seen that ${N^{\text{max}}_{\text{MC,prim}} < N^{\text{max}}_{\text{prim}}}$ and ${N^{\text{max}}_{\text{MC,dual}} < N^{\text{max}}_{\text{dual}}}$ for some of the examples, highlighted in red in Table \ref{tab:res}. This either means that the certification method is conservative or that the MC simulations are optimistic, (or both). However, since the certification method provides a region in parameter space for which the worst-case number of iterations is obtained, a parameter in the worst-case region for each example was extracted and by applying Algorithm \ref{alg:primal-as} to the resulting QP it could be proven that the certification method did not provide a conservative result. Instead, the discrepancies are due to MC simulations not being able to cover the parameter space densely enough with samples during the allotted time. Even if more samples could be taken to improve the MC results, this would require more time than the certification method and, still, there are no guarantees for sufficient coverage for any finite number of samples. This underlines an important advantage of the proposed certification method compared to MC simulations, namely that the proposed method covers a continuum of points, which becomes increasingly beneficial as the dimension of the parameter space increases. 

\begin{remark}
  The execution time $t^{\text{cert}}$ is based on a implementation of Algorithm \ref{alg:main} in MATLAB. Modifications to the implementation, such as low-rank modifications and parallelizing computations, are expected to significantly reduce $t^{\text{cert}}$. 
\end{remark}
\subsection{Affine approximations of quadratic inequalities}
The affine outer-approximations of quadratic constraints, described in Section \ref{ssec:outer-approx}, were tested by using Algorithm \ref{alg:main} with and without these relaxations on the problems which lead to quadratic partitioning of parameter space, i.e.,  the contrived mpQP, DC motor and ATFI-16 aircraft example. Table \ref{tab:relax} summarizes the result, where it can be seen that approximating the quadratic constraints results in the final partition containing more regions, given by $N^{\text{reg}}$, for all of the examples. This is expected since, as is discussed in Section \ref{ssec:outer-approx}, the relaxations might lead to redundant regions. For the contrived mpQP, the relaxation results in an upper bound on the number of QP iterations $N^{\text{max}}_{\text{iter}}$ of 6 instead of the tight upper bound 4. However, for both the DC motor and ATFI-F16 example the upper bounds provided by the relaxation coincide with the tight upper bound. Table \ref{tab:relax} also shows that, for large problems, the computation time $t^{\text{cert}}$ for the certification can be reduced significantly by forming affine outer-approximations of the quadratic constraints. In conclusion, relaxing quadratic constraints with the method described in Section \ref{ssec:outer-approx} can provide good, even tight, upper bounds on worst-case behaviour while reducing the certification time for large problems. 

\begin{table}[h]
  \centering
  \begin{tabular}{c || c| c| c}
  & $N^{\text{max}}_{\text{iter}}$ & $N^{\text{reg}}$ & $t^{\text{cert}}[s]$ \\ \hline 
  Contrived mpQP& $6/4$ &$15/6$& $0.76/0.08 $\\ %\hline 
  DC motor & $14/14$& $2309/1765$& $46/263$\\ %\hline 
  ATFI-F16 & $21/21$&$41971/31831$ & $541/4689$ 
  \end{tabular}
  \caption{Comparison of the certification method when linear outer-approximations of quadratic constraints are/aren't used.}
  \label{tab:relax}
\end{table}

  \section{CONCLUSION AND FUTURE WORK}
  In this paper we have presented a method which extends, and unifies, complexity certification results for active-set QP and LP methods. 
The method computes \textit{exactly} which sequence of working-set changes, as a function of the parameters in an mpQP, a primal active-set QP algorithm will undergo to find an optimum. This can be used to determine an upper bound on the number of QP iterations the algorithm will need when it is applied online,  which is of importance in the context of real-time MPC where hard real-time requirements have to be fulfilled. The method partitions the parameter space into regions, defined by affine and quadratic inequalities, representing parameter sets which generate the same sequence of working-set changes to reach a solution.  
Furthermore, by considering positive semi-definite QPs, the proposed method poses previous complexity certification results for primal and dual active-set QP methods, as well as active-set LP methods, in a unified framework. 
The proposed method was successfully applied to a set of linear MPC problems to illustrate how it can be used to determine the worst-case number of iterations needed by a primal and a dual active-set algorithm online.

Future work includes using the framework to compare the worst-case number of FLOPs different active-set algorithms result in, e.g., the difference between different range-space and null-space methods.  

\bibliographystyle{IEEEtran}
\bibliography{lib.bib}
\end{document}